\documentclass[12pt]{amsart}
\usepackage{setspace}
\usepackage[english]{babel}
\usepackage{bbold}
\usepackage{stmaryrd}
\usepackage{verbatim}
\usepackage[utf8x]{inputenc}
\usepackage{amsmath}
\numberwithin{equation}{subsection}
\usepackage{amssymb}
\usepackage{amsthm}
\usepackage{graphicx}
\usepackage{mathtools}
\usepackage[T1]{fontenc}
\usepackage{amsthm}
\usepackage{setspace}
\usepackage{amsfonts}
\usepackage{eucal}
\usepackage{mathrsfs}
\usepackage{textcomp}
\usepackage{pictexwd,dcpic}
\usepackage{tikz}
\usepackage{tqft}
\usepackage{pigpen}
\usepackage{xcolor}
\usepackage{pdfpages}
\usepackage{hyperref}
\newtheorem*{theorem-non}{Theorem}
\newtheorem*{theorem-nonf}{Théorème}
\newtheorem*{corollary-non}{Corollary}
\newtheorem*{corollary-nonf}{Corollaire}

\newtheorem{lmm}[subsection]{Lemma}
\newtheorem{thm}[subsection]{Theorem}
\newtheorem{prop}[subsection]{Proposition}
\newtheorem{cor}[subsection]{Corollary}

\theoremstyle{definition}
\newtheorem{defn}[subsection]{Definition}
\newtheorem{context}[subsection]{Context}
\newtheorem{ex}[subsection]{Example}
\newtheorem*{ex-non}{Example}
\newtheorem{rmk}[subsection]{Remark}

\newtheorem*{rmk-non}{Remark}
\newtheorem*{rmk-nonf}{Remarque}

\newtheorem{notation}[subsection]{Notation}

\title{}
\date{ }

\newcommand{\drl}[0]{\llparenthesis}
\newcommand{\drr}[0]{\rrparenthesis}
\newcommand{\dsl}[0]{\llbracket}
\newcommand{\dsr}[0]{\rrbracket}
\newcommand{\ital}[1]{\textit{#1}}

\newcommand{\red}[1]{\textcolor{red}{#1}}

\newcommand{\Bc}[0]{\mathcal{B}}
\newcommand{\Cc}[0]{\mathcal{C}}
\newcommand{\Ec}[0]{\mathcal{E}}
\newcommand{\Fc}[0]{\mathcal{F}}

\newcommand{\Lc}[0]{\mathcal{L}}
\newcommand{\Mcc}[0]{\mathcal{M}}

\newcommand{\Oc}[0]{\mathcal{O}}
\newcommand{\Pc}[0]{\mathcal{P}}

\newcommand{\Uc}[0]{\mathcal{U}}
\newcommand{\Vc}[0]{\mathcal{V}}

\newcommand{\Es}[0]{\mathscr{E}}
\newcommand{\Is}[0]{\mathscr{I}}
\newcommand{\Ns}[0]{\mathscr{N}}
\newcommand{\C}[0]{\mathbb{C}}

\newcommand{\Pbb}[0]{\mathbb{P}}
\newcommand{\Q}[0]{\mathbb{Q}}
\newcommand{\V}[0]{\mathbb{V}}
\newcommand{\Z}[0]{\mathbb{Z}}

\newcommand{\Hu}[0]{\textup{H}}
\newcommand{\Vu}[0]{\textup{V}}

\newcommand{\oo}{\infty}

\newcommand{\Rlv}[0]{\mathcal{R}^{\ell,\vee}}

\newcommand{\Mv}[0]{\mathcal{M}^{\vee}}

\newcommand{\Ql}[1]{\mathbb{Q}_{\ell#1}}
\newcommand{\aff}[2]{\mathbb{A}^{#1}_{#2}}
\newcommand{\proj}[2]{\mathbb{P}^{#1}_{#2}}

\newcommand{\ag}[1]{[\mathbb{A}^1_{#1}/\mathbb{G}_{m,#1}]}

\newcommand{\bu}[0]{\textup{B}\mathbb{U}}

\newcommand{\shvl}[0]{\textup{\textbf{Shv}}_{\mathbb{Q}_{\ell}}}

\newcommand{\md}[0]{\textup{\textbf{Mod}}}

\newcommand{\tlgmprojm}[1]{\textup{LG}_{(\mathbb{P}(#1),\mathcal{O}(1))}^{\boxplus}}
\newcommand{\tlgmprojmop}[1]{\textup{LG}_{(\mathbb{P}(#1),\mathcal{O}(1))}^{\boxplus,\textup{op}}}

\newcommand{\tlgmrm}[1]{\textup{LG}_{(#1,\Ec_{#1})}^{\boxplus}}
\newcommand{\tlgmrmop}[1]{\textup{LG}_{(#1,\Ec_{#1})}^{\boxplus,\textup{op}}}
\newcommand{\tlgmp}[1]{\textup{LG}_{(#1,\Ec^+_{#1})}}
\newcommand{\tlgmpm}[1]{\textup{LG}_{(#1,\Ec^+_{#1})}^{\boxplus}}
\newcommand{\tlgmpop}[1]{\textup{LG}_{(#1,\Ec^+_{#1})}^{\textup{op}}}
\newcommand{\tlgmpmop}[1]{\textup{LG}_{(#1,\Ec^+_{#1})}^{\boxplus,\textup{op}}}

\newcommand{\dgcatm}[1]{\textup{\textbf{dgCat}}^{\textup{idm}}_{#1}}

\newcommand{\dgcatmo}[1]{\textup{\textbf{dgCat}}^{\textup{idm},\otimes}_{#1}}

\newcommand{\cohbsupp}[2]{\textup{\textbf{Coh}}^{b}_{#2}(#1)}
\newcommand{\cohb}[1]{\textup{\textbf{Coh}}^{b}(#1)}
\newcommand{\perf}[1]{\textup{\textbf{Perf}}(#1)}
\newcommand{\perfsupp}[2]{\textup{\textbf{Perf}}_{#2}(#1)}
\newcommand{\rperf}[2]{\textup{\textlangle}#2_*\textup{\textbf{Perf}}(#1)\textup{\textrangle}}
\newcommand{\qcoh}[1]{\textup{\textbf{QCoh}}(#1)}

\newcommand{\cohbp}[2]{\textup{\textbf{Coh}}^b(#1)_{\textup{\textbf{Perf}}(#2)}}

\newcommand{\rcoh}[2]{\textup{\textlangle} #2_* \textup{\textbf{Coh}}^{b}(#1) \textup{\textrangle}}
\newcommand{\cohs}[1]{\textup{Coh}^{s}(#1)}
\newcommand{\sing}[1]{\textup{\textbf{Sing}}(#1)}
\newcommand{\singsupp}[2]{\textup{\textbf{Sing}}_{#2}(#1)}
\newcommand{\singc}[1]{\textup{\textbf{Sing}}^{\textup{coh}}(#1)}

\newcommand{\ising}[2]{\textup{\textlangle} #2^* \textup{\textbf{Sing}}(#1) \textup{\textrangle}}
\newcommand{\rsing}[2]{\textup{\textlangle} #2_* \textup{\textbf{Sing}}(#1) \textup{\textrangle}}
\newcommand{\mf}[2]{\textup{\textbf{MF}}(#1,#2)}
\newcommand{\mfc}[2]{\textup{\textbf{MF}}^{\textup{coh}}(#1,#2)}
\newcommand{\singcf}[1]{\mathscr{S}^{\textup{coh}}_{(#1,\Ec_{#1}^+)}}

\newcommand{\singpf}[1]{\textup{\textbf{Sing}}_{(#1,\Ec_{#1}^+)}}
\newcommand{\singpfm}[1]{\textup{\textbf{Sing}}_{(#1,\Ec_{#1}^+)}^{\otimes}}
\newcommand{\singf}[1]{\textup{\textbf{Sing}}_{(#1,\Ec_{#1})}}
\newcommand{\singfm}[1]{\textup{\textbf{Sing}}_{(#1,\Ec_{#1})}^{\otimes}}
\newcommand{\singprojf}[1]{\textup{\textbf{Sing}}_{(\mathbb{P}(#1),\Oc(1))}}
\newcommand{\singprojfm}[1]{\textup{\textbf{Sing}}_{(\mathbb{P}(#1),\Oc(1))}^{\otimes}}
\newcommand{\singprojcf}[1]{\mathscr{S}^{\textup{coh}}_{(\mathbb{P}(#1),\Oc(1))}}

\newcommand{\hh}[1]{\textup{\textbf{HH}}_{\bullet}(#1)}
\newcommand{\hp}[1]{\textup{\textbf{HP}}_{\bullet}(#1)}
\newcommand{\hoch}[2]{(\textup{Hoch}(#1),#2)}

\newcommand{\rh}[2]{\widehat{\textup{RH}}^{-1}(#1,#2)}

\title{\parbox{\linewidth}{\centering On some (co)homological invariants of coherent matrix factorizations}}
\date{}
\author{Massimo Pippi}
\email{m.pippi@ucl.ac.uk}
\address{Department of Mathematics, University College London, Gordon Street 25, London WC1H 0AY, United Kingdom}

\begin{document}
\maketitle
\begin{abstract}
We provide an equivalence between the dg category of coherent matrix factorizations and a certain dg category of absolute singularities. As an application, we compute the $\ell$-adic cohomology of the dg category of coherent matrix factorizations, as well as its Hochschild and periodic cyclic homologies (these last two only in the affine case).
\end{abstract}
\tableofcontents

\section{Introduction} 
Matrix factorizations are non-commutative spaces (a.k.a. dg categories) associated to varieties with a potential $h:Y\rightarrow \aff{1}{}$. They were first introduced by Eisenbud in \cite{ei80} in order to study hypersurface singularities. Buchweitz (\cite{buch87}) and Orlov (\cite{orl04}) independently found out that these are related to the dg category of singularities of the fiber over $0$ of $h$. This is defined as the Drinfeld/Verdier quotient
\begin{equation*}
    \sing{h^{-1}(0)}=\frac{\cohb{h^{-1}(0)}}{\perf{h^{-1}(0)}}.
\end{equation*}
In particular, they showed that these two dg categories are equivalent if $Y$ and $h$ are regular.
In other words, matrix factorizations are related to those complexes with coherent cohomology over $h^{-1}(0)$ that are not quasi-isomorphic to a perfect one.

This connection can be explained by means of the main theorem in \cite{ei80}. Assume that $Y=Spec(Q)$ is a regular local ring of finite Krull dimension. Let $X=Spec(Q/h)$. A classical result due to Auslander-Buchsbaum/Serre (see \cite{ab56} and \cite{se55}) asserts that $X$ is regular if and only if it has finite global dimension, i.e. any finitely generated module admits a finite resolution by projective modules of finite type. Eisenbud proved that if $X$ fails to be regular, the resolutions of finitely generated $Q/h$ modules eventually become $2$ periodic. This $2$ periodic part of the complex is essentially the datum of two projective $Q$-modules of finite type, together with two $Q$-linear morphisms from one another whose compositions are equal to multiplication by $h$. That is, the two periodic part of the complex is essentially the datum of a matrix factorization.

Various generalizations of the dg category of matrix factorizations and of dg categories of singularities have been introduced in the literature.

The main theme of this note deals with dg categories of matrix factorizations defined for pairs $(Y,h)$ where $Y$ is allowed to be singular.

When $Y$ is a regular scheme, there is only one candidate for what should be the category of matrix factorizations of $h:Y\rightarrow \aff{1}{}$ and it coincides with the dg category of absolute singularities of $h^{-1}(0)$. However, as soon as we allow $Y$ to be singular, (at least) two non equivalent definitions are possible. The first possibility is to consider the dg category of relative singularities of $i:h^{-1}(0)\rightarrow Y$ (\cite{bw12}, \cite{ep15}). The pushforward induces a well defined dg functor
\begin{equation*}
    i_*:\sing{h^{-1}(0)}\rightarrow \sing{Y}
\end{equation*}
and the dg category of relative singularities is defined as its kernel:
\begin{equation*}
    \sing{Y,h}:=Ker \bigl ( i_*:\sing{W^{-1}(0)}\rightarrow \sing{Y}\bigr ).
\end{equation*}
\begin{comment}
\begin{equation*}
    \sing{Z\subseteq Y}=Ker(i_*:\sing{Z}\rightarrow \sing{Y}).
\end{equation*}
As far as the author is aware, it was first introduced by Burke and Walker in \cite{bw12}. This dg category can alternatively be described as the dg quotient
\begin{equation*}
    \sing{Z\subseteq Y}\simeq \cohbp{Z}{Y}/\perf{Z},
\end{equation*}
where $\cohbp{Z}{Y}$ is the sub dg category of $\cohb{Z}$ spanned by those objects $\Bc$ such that $i_*\Bc$ lies in $\perf{Y}$. This is the version that was considered in the previous section.
\end{comment}
The pullback $i^*$ also induces a dg functor between the dg categories of singularities of $Y$ and $h^{-1}(0)$. Therefore, as a second possibility, we can consider the dg quotient
\begin{equation*}
    \singc{Y,h}= Coker(i^*:\sing{Y}\rightarrow \sing{h^{-1}(0)}),
\end{equation*}
also known as the dg category of coherent relative singularities. As far as the author is aware, this dg category was first considered by Efimov and Positselski in \cite{ep15}.

These two relative versions of the dg category of singularities also admit a description in terms of matrix factorizations (\cite{bw12}, \cite{ep15}, \cite{brtv18}, \cite{p19}). Roughly, $\sing{Y,h}$ is equivalent to the dg category $\mf{Y}{h}$ whose objects are quadruples $(E_0,E_1,\delta_0:E_0\rightarrow E_1,\delta_1:E_1\rightarrow E_0)$ where $E_0$ and $E_1$ are $Y$-vector bundles of finite rank and $\delta_0\circ \delta_1=h=\delta_1\circ \delta_0$, while $\singc{Y,h}$ is equivalent to the dg category $\mfc{Y}{h}$ whose objects are quadruples $(F_0,F_1,\partial_0:F_0\rightarrow F_1,\partial_1:F_1\rightarrow F_0)$ where $F_0$ and $F_1$ are coherent $\Oc_Y$ modules and $\partial_0\circ \partial_1=h=\partial_1\circ \partial_0$. This also explains the choice of terminology for $\singc{Y,h}$.

\begin{rmk}
It is immediate to see that $\sing{Y,h}\simeq \singc{Y,h}\simeq \sing{h^{-1}(0)}$ when $Y$ is regular, as in this case $\sing{Y}\simeq 0$.
\end{rmk}

An interesting feature of non commutative spaces is that, while they do not carry any topology, it is still possible to define interesting (co)homological invariants. 

The main purpose of this note is the computation of some (co)homological invariants of the dg category of coherent matrix factorizations. As far as the author is aware, these (co)homological invariants were previously computed only in the case where the ambient space is regular (\cite{dy11}, \cite{pr11}, \cite{lp13},  \cite{se13}, \cite{ef18}, \cite{brtv18}, \cite{p20}).

The main observation is that sometimes it is possible to identify the category of coherent matrix factorizations $\mfc{Y}{h}$ with the category of (usual, i.e. locally free) matrix factorizations $\mf{\Vc}{W}$ of a pair $(\Vc,W)$ where $\Vc$ is regular. This is done by studying certain (well known) cohomological operators on the two sides of a reduction of codimension equivalence for matrix factorizations (see \cite{or06}, \cite{bw15}).

More precisely, let $X$ be a quasi compact and separated scheme of finite type over a regular local noetherian ring $S$ of finite Krull dimension. Let $\Ec_X$ be a vector bundle over $X$ and let $\Lc_X$ be a line bundle over $X$. Moreover, let $p^+=(p,s)$ be a global section of $\Ec_X^{+}=\Ec_X\oplus \Lc_X$.
Recall that, for any locally noetherian (derived) scheme $Z$ (such that $\Oc_Z \in \cohb{Z}$), we can define its dg category of absolute singularities:
\begin{equation*}
    \sing{Z}:=\frac{\cohb{Z}}{\perf{Z}}.
\end{equation*}
Starting from $(X,p^+)$, we can consider the following two dg categories of relative singularities
\begin{equation*}
    \sing{X,p^+}:=Ker\bigl ( f_*:\sing{\Vu(p^+)}\rightarrow \sing{X}\bigr ),
\end{equation*}
\begin{equation*}
    \sing{X,p}:=Ker\bigl ( g_*:\sing{\Vu(p)}\rightarrow \sing{X}\bigr ),
\end{equation*}
where $f:\Vu(p^+)\rightarrow X$ and $g:\Vu(p)\rightarrow X$ are the lci closed embedding determined by $p^+$ and $p$ respectively. 
A theorem of Orlov, refined by Burke and Walker, tells us that we have equivalences
\begin{equation*}
    \sing{X,p^+}\simeq \sing{\Pbb(\Ec_X^{+,\vee}),W_{p^+}}, \hspace{0.5cm}
    \sing{X,p}\simeq \sing{\Pbb(\Ec_X^{\vee}),W_{p}}.
\end{equation*}
Here $\Pbb(\Ec_X^{+,\vee})=Proj_X(Sym_{\Oc_X}(\Ec_X^+))$ (resp. $\Pbb(\Ec_X^{\vee})=Proj_X(Sym_{\Oc_X}(\Ec_X))$), $W_{p^+}$ (resp. $W_p$) denotes the global section of $\Oc_{\Pbb(\Ec_X^{+,\vee})}(1)$ (resp. $\Oc_{\Pbb(\Ec_X^{\vee})}(1)$) determined by $p^+$ (resp. $p$) and
\begin{equation*}
    \sing{\Pbb(\Ec_X^{+,\vee}),W_{p^+}}:=Ker\bigl ( u_*:\sing{\Vu(W_{p^+})}\rightarrow \sing{\Pbb(\Ec_X^{+,\vee})}\bigr ),
\end{equation*}
\begin{equation*}
    \sing{\Pbb(\Ec_X^{\vee}),W_{p}}:=Ker\bigl ( v_*:\sing{\Vu(W_{p})}\rightarrow \sing{\Pbb(\Ec_X^{\vee})}\bigr ),
\end{equation*}
where $u:\Vu(W_{p^+})\rightarrow \Pbb(\Ec_X^{+,\vee})$, $v:\Vu(W_{p})\rightarrow \Pbb(\Ec_X^{\vee})$.
There are lci closed embeddings
\begin{equation*}
    i:\Vu(p^+)\hookrightarrow \Vu(p), \hspace{0.5cm} r: \Vu(W_{p})\hookrightarrow \Vu(W_{p^+})
\end{equation*}
that induce adjunctions
\begin{equation*}
    i^*:\sing{X,p}\rightleftarrows \sing{X,p^+}:i_*, \hspace{0.5cm} r_*:\sing{\Pbb(\Ec_X^{\vee}),W_{p}}\rightleftarrows \sing{\Pbb(\Ec_X^{+,\vee}),W_{p^+}}:r^!.
\end{equation*}
The crucial observation of this note is that these adjunctions correspond (up to a twisted shift):
\begin{theorem-non}{\textup{\textbf{\ref{equivalence adjunction}}}}
The adjunctions 
\begin{equation*}
 i^*\bigl ((-)\otimes \Lc_{\Vu(p)}  \bigr )[-1] :\sing{X,p}\rightarrow \sing{X,p^+}: \bigl (i_*(-)\otimes \Lc_{\Vu(p)}^{\vee}\bigr )[1]
\end{equation*}
\begin{equation*}
     r_*:\sing{\Pbb(\Ec_X^{\vee}),W_p}\rightarrow \sing{\Pbb(\Ec_X^{+,\vee}),W_{p^+}}:r^!
\end{equation*}  
identify under the equivalences above.
\end{theorem-non}
The counit of $(i^*,i_*)$ (resp. $(r_*,r^!)$) induces a \textit{cohomological operator} $\varepsilon$ (resp. $\chi$) on $\sing{X,p^+}$ (resp. $\sing{\Pbb(\Ec_X^{+,\vee}),W_{p^+}}$). See Definitions \ref{defn E-operator} and \ref{defn BW operator}. As an immediate consequence of the Theorem above, these operators agree under the equivalences above.

We will next study the torsion objects of these operators, i.e. objects $\Pc$ in $\sing{X,p^+}$ (resp. $\sing{\Pbb(\Ec_X^{+,\vee}),W_{p^+}}$) such that $\varepsilon^n\sim 0$ (resp. $\chi^n\sim 0$) for some $n>>0$. Clearly, these agree under the equivalence $\sing{X,p^+}\simeq \sing{\Pbb(\Ec_X^{+,\vee}),W_{p^+}}$.

We will show that (see Proposition \ref{E-torsion = image pullback})
\begin{equation*}
    \sing{X,p^+}^{\varepsilon-tors}\simeq \ising{X,p}{i},
\end{equation*}
where $\ising{X,p}{i}$ is the thick sub category of $\sing{X,p^+}$ generated by the image of the pullback $i^*:\sing{X,p}\rightarrow \sing{X,p^+}$,
and that (see Proposition \ref{BW-torsion = image r_*})
\begin{equation*}
    \sing{\Pbb(\Ec_X^{+,\vee}),W_{p^+}}^{\chi-tors}\simeq \rsing{\Pbb(\Ec_X^{\vee}),W_p}{r},
\end{equation*}
where $\rsing{\Pbb(\Ec_X^{\vee}),W_p}{r}$ denotes the thick sub category generated by the image of the pushforward $r_*:\sing{\Pbb(\Ec_X^{\vee}),W_p}\rightarrow \sing{\Pbb(\Ec_X^{+,\vee}),W_{p^+}}$. 

Furthermore, when $X$ and $p^+$ are regular, we will identify $\rsing{\Pbb(\Ec_X^{\vee}),W_p}{r}$ with the dg category $\singsupp{\Pbb(\Ec_X^{+,\vee}),W_{p^+}}{\Vu(W_p)}\simeq \singsupp{\Vu(W_{p^+})}{\Vu(W_p)}$, which is the full sub category of $\sing{\Pbb(\Ec_X^{+,\vee}),W_{p^+}}\simeq \sing{\Vu(W_{p^+})}$ spanned by those objects whose support is contained in $\Vu(W_p)$.

Therefore, assuming that $X$ and $p^+$ are regular and putting all pieces together, we get 

\begin{corollary-non}\textup{\textbf{\ref{main equivalence}}}
The equivalence $\sing{X,p^+}\simeq \sing{\Pbb(\Ec_X^{+,\vee}),W_{p^+}}$ induces an equivalence of exact sequences in $\dgcatm{S}$:
\begin{equation*}
    \begindc{\commdiag}[20]
    \obj(-60,10)[1]{$\ising{\Vu(p)}{i}$}
    \obj(0,10)[2]{$\sing{\Vu(p^+)}$}
    \obj(60,10)[3]{$\singc{\Vu(p),s}$}
    \obj(-60,-10)[4]{$\singsupp{\Vu(W_{p^+})}{\Vu(W_p)}$}
    \obj(0,-10)[5]{$\sing{\Vu(W_{p^+})}$}
    \obj(60,-10)[6]{$\sing{\Vu(W_{p^+})-\Vu(W_p)}.$}
    \mor{1}{2}{$$}[\atright,\injectionarrow]
    \mor{2}{3}{$$}
    \mor{1}{4}{$$}
    \mor{2}{5}{$$}
    \mor{3}{6}{$$}
    \mor{4}{5}{$$}[\atright,\injectionarrow]
    \mor{5}{6}{$$}
    \enddc
\end{equation*}

In particular,
\begin{equation*}
    \singc{\Vu(p),s}=\frac{\sing{\Vu(p^+)}}{\ising{\Vu(p)}{i}}\simeq \frac{\sing{\Vu(W_{p^+})}}{\rsing{\Vu(W_{p})}{r}}\simeq \sing{\Vc,W_{p^+|\Vc}},
\end{equation*}
where $\Vc=\Pbb(\Ec_X^+)-\Pbb(\Ec_X)$ is a regular scheme.
\end{corollary-non}
The dg category on the left hand side is the dg category of coherent relative singularities of $\Vu(p^+)\hookrightarrow \Vu(p)$ . The main advantage of the equivalence above is the identification of $\singc{\Vu(p),s}$ with the category of matrix factorizations of a pair $(\Vc,W_{p^+|\Vc})$ with $\Vc$ regular.
\begin{rmk}
Let $h:Y\rightarrow \aff{1}{}$ be a variety with a potential. If $Y$ is regular, then $\mf{Y}{h}$ has nice properties, e.g. it satisfies Thomason's localization. When $Y$ is not regular this is no longer true (see \cite{ep15} and \cite{ef18}): Thomason's localization fails in general, periodic cyclic homology may be infinite dimensional $\dots$ It is true however that $\mfc{Y}{h}$ has nice properties. The equivalence above might shed some light on different behaviours of the two versions of MF: being equivalent to a dg category of matrix factorizations for a pair where the ambient space is regular, the coherent version inherits all its nice properties.
\end{rmk}

\begin{rmk}\label{mfcoh vanishies for almost all c's}
Assume that $g:X=Spec(A)\rightarrow \aff{1}{\C}$ is a variety with a potential. It is well known (see e.g. \cite{orl04}) that, if $X$ is regular, $\mf{X}{g}$ is non zero if and only if $g^{-1}(0)$ is singular. In particular, only finitely many $\mf{X}{g-c}$, $c\in \C$, do not vanish in this case. If $X$ is not supposed to be regular, the author was not aware of any argument guaranteeing that $\mfc{X}{g-c}$ is non zero only for finitely many $c\in \C$. This can be seen by means of the previous equivalence. Assume that $(f_1,\dots,f_m,g)\in \C[x_1,\dots,x_n]$ is a regular sequence and let $A=\frac{\C[x_1,\dots,x_n]}{(f_1,\dots,f_m)}$. Since $X$ may be singular, it is not true that only finitely many values of $g:X\rightarrow \aff{1}{\C}$ are critical. Nevertheless, it is not true that $\mfc{X}{g-c}$ is non zero if $g^{-1}(c)$ is singular neither. The equivalence above translates in this context as
\begin{equation*}
    \mfc{X}{g}\simeq \mf{\aff{n+m}{\C}}{W},
\end{equation*}
where $W=f_1\cdot x_{n+1}+\dots+f_m\cdot x_{n+m}+g$. Also notice that, for any $c\in \C$, the assignment $g\rightsquigarrow W$ is compatible with translations, i.e. $g-c\rightsquigarrow W-c$. More generally, we thus have that
\begin{equation*}
    \mfc{X}{g-c}\simeq \mf{\aff{n+m}{\C}}{W-c}.
\end{equation*}
As the dg category on the right hand side is non zero only for finitely many values of $c\in \C$, the same holds true for the dg category of coherent matrix factorizations.

This phenomenon can be explained as, roughly, $\mfc{X}{g-c}$ is non zero if and only if the points of the fiber $g^{-1}(c)\subseteq Y$ are "more singular" when considered inside $g^{-1}(c)$ rather then when considered in $Y$. See also Definition \ref{relative critical locus}.
\end{rmk}
\begin{ex}
We will try to explain the phenomenon of the previous remark via an example. Let $g=x_2^2+x_3^2+x_4^2+x_5^2, f=x_1\cdot x_2+x_3\cdot x_4 \in \C[x_1,\dots,x_5]$ and $X=\Vu(f)\subseteq \aff{5}{\C}$. Notice that all points $(0,0,0,0,\pm \sqrt{c}) \in g^{-1}(c)$ are singular. In particular, all fibers $g^{-1}(c)$ are singular.

In this example,
\begin{equation*}
    W=(x_1\cdot x_2 + x_3\cdot x_4)\cdot x_6 +x_2^2+x_3^2+x_4^2+x_5^2:\aff{6}{\C}\rightarrow \aff{1}{\C}.
\end{equation*}
and one easily sees that $0$ is the only critical value of $W$. In particular,
\begin{equation*}
    \mfc{X}{g-c}\simeq \mf{\aff{6}{\C}}{W-c}\simeq 0
\end{equation*}
for all $c\neq 0$.

Alternatively, one can argue that the reason for which $\mfc{X}{g-c}\simeq 0$ if $c\neq 0$ is that $0$ is the only relative critical value of $g:X\rightarrow \aff{1}{\C}$ (in the sense of Definition \ref{relative critical locus}).
\end{ex}
\begin{rmk}
Let $g:X\rightarrow \aff{1}{\C}$ be as in Remark \ref{mfcoh vanishies for almost all c's}.
Notice that also $\mf{X}{g-c}$, the dg category locally free matrix factorizations, vanishes  for almost all $c\in \C$. This is an immediate consequence of Remark \ref{mfcoh vanishies for almost all c's} and of the existence of a fully faithful dg functor (see \cite{ep15})
\begin{equation*}
    \mf{X}{g-c}\subseteq \mfc{X}{g-c}.
\end{equation*}
\end{rmk}

In the last section we apply the equivalence above to compute certain (co)homological invariants of $\singc{\Vu(p),s}$, namely its $\ell$-adic cohomology (as defined in \cite{brtv18}), its Hochschild homology and its periodic cyclic homology (in the affine case). These results can be seen as attempts to generalize previous established theorems: for the $\ell$-adic cohomology, the main theorems in \cite{brtv18} and in \cite{p20}; for the Hochschild and periodic cyclic homologies, theorems in \cite{dy11} \cite{ef18}, \cite{lp13}, \cite{pr11}, \cite{se13}.
\begin{rmk}
For more details on our notation, we invite the reader to consult section \S \ref{applications}.
\end{rmk}
\begin{theorem-non}{\textup{\textbf{\ref{l-adic cohomology mf coh}}}}
Let $(X,p^+)$ be an LG model over $(Y,\Ec_Y^+)$ such that $X$ is a regular scheme and $p^+$ is a regular section. Then
\begin{equation*}
   \Rlv_{Y}(\singc{\Vu(p),s})\simeq \Rlv_{Y}(\mfc{\Vu(p)}{s}\simeq q_*\Phi_{(\Vc,W_{p^+|\Vc})}^{\textup{mi}}(\Ql{}(\beta))[-1],
\end{equation*}
where $q:\Uc=\Vu(W_{p^+|\Vc}) \rightarrow Y$ is the canonical morphism, $\Rlv_{Y}$ denotes the "$\ell$-adic cohomology of dg categories" $\oo$-functor and $\Phi_{(\Vc,W_{p^+|\Vc})}^{\textup{mi}}(\Ql{}(\beta))$ denotes the $\ell$-adic sheaf of monodromy invariant vanishing cycles (see \cite{p20}).
\end{theorem-non}

\begin{theorem-non}{\textup{\textbf{\ref{computation HH and HP mf coh}}}}
Let $f_1,\dots,f_m,g \in \C[x_1,\dots,x_m]$ be a regular sequence. Let $X=Spec\bigl (\frac{\C[x_1,\dots,x_n]}{(f_1,\dots,f_m)}\bigr )$ and let $W=g+\sum_{k=1}^mf_k\cdot x_{n+k}\in \C[x_1,\dots,x_{n+m}]$.
There is an equivalence of $\Z/2\Z$-graded bundles with connections on $Spf(\C \drl u \drr)$
\begin{equation*}
    \bigoplus_{c\in \C} \hp{\mfc{X}{g-c},\nabla_u^{GM}}\simeq \bigl ( \Hu^{\bullet}_{Zar}(\Omega^{\bullet}_{\aff{n+m}{\C}}\drl u \drr, -dW+ud_{dR}), \nabla_u^{dR}=\frac{d}{du}+\frac{\Gamma}{u}+\frac{W}{u^2} \bigr ),
\end{equation*}
where $\Gamma_{|\Omega^p_{\aff{n+m}{\C}}}=-\frac{p}{2}\cdot id$, and an equivalence of $\Z/2\Z$-graded vector spaces
\begin{equation*}
    \bigoplus_{c\in \C}\hh{\mfc{X}{g-c}}\simeq \Hu^{\bullet}_{Zar}(\aff{n+m}{\C},(\Omega^{\bullet}_{\aff{n+m}{\C}},-dW)).
\end{equation*}
Moreover, if the relative critical locus of $g$ (see Definition \ref{relative critical locus}) is contained in the fiber over $0$, then
\begin{equation*}
   \hp{\mfc{X}{g},\nabla_u^{GM}}\simeq \bigl ( \Hu^{\bullet}_{Zar}(\Omega^{\bullet}_{\aff{n+m}{\C}}\drl u \drr, -dW+ud_{dR}), \nabla_u^{dR}=\frac{d}{du}+\frac{\Gamma}{u}+\frac{W}{u^2} \bigr ),
\end{equation*}
and 
\begin{equation*}
    \hh{\mfc{X}{g}}\simeq \Hu^{\bullet}_{Zar}(\aff{n+m}{\C},(\Omega^{\bullet}_{\aff{n+m}{\C}},-dW)).
\end{equation*}
\end{theorem-non}
Combining the theorem above with the algebraic computation of vanishing cycles (\cite{sa12}, \cite{sasa14}), we get
\begin{theorem-non}{\textup{\textbf{\ref{HP and vanishing cycles}}}}
With the same notation as above, there is an equivalence of $\Z/2\Z$-graded vector bundles with connections of the punctured disk, 
\begin{equation*}
    \bigoplus_{c\in \C}\hp{\mfc{X}{g-c},\nabla_u^{GM}}\simeq \bigoplus_{c\in \C}\widehat{\Es}^{-\frac{c}{u}}\otimes_{\C \drl u \drr} \rh{\Hu^{\bullet-1}_{an}(W^{-1}(c),\Phi_{W-c}(\C_{\aff{n+m}{\C}}))}{T\cdot (-1)^{\bullet}}.
\end{equation*}
Moreover, if the relative critical locus of $g$ is contained in the fiber over $0$, then
\begin{equation*}
    \hp{\mfc{X}{g},\nabla_u^{GM}}\simeq  \rh{\Hu^{\bullet-1}_{an}(W^{-1}(0),\Phi_{W}(\C_{\aff{n+m}{\C}}))}{T\cdot (-1)^{\bullet}}.
\end{equation*}
\end{theorem-non}
\begin{rmk}
Matrix factorizations play a crucial role in
Toën-Vezzosi's approach to the Bloch's conductor formula conjecture (see \cite{b85} and \cite{tv19}). In \cite{tv17} the authors suggest that applying their machine to $\mfc{Z}{\pi_Z}$ (for a suitable pair $(Z,\pi_Z)$) might be useful to remove the unipotency hypothesis in their categorical version of the Bloch's conductor formula.

Notice however that Toën and Vezzosi are currently investigating a second approach to the Bloch's conductor conjecture which does not make use of coherent matrix factorizations. See \cite{v19}.
\end{rmk}
\begin{center} 
   \ital{Conventions and notations}
\end{center}
\begin{itemize}
    \item We will use the theory of $\oo$-categories. See \cite{lu09} and \cite{lu17}.
    \item We will use interchangeably the terminology "matrix factorizations" and "dg category of (relative) singularities". Indeed, it is known after Buchweitz (\cite{buch87}), Orlov (\cite{orl04}), Burke-Walker (\cite{bw12},\cite{bw15}) and Efimov-Positselski (\cite{ep15}) (see also \cite{p19}) that these coincide. In particular, for us "matrix factorizations" will be what are sometimes called "locally free matrix factorizations" in the literature.
    \item Similarly, we will say "coherent matrix factorizations" and "category of coherent relative singularities" to identify the same category. This is possible after \cite{ep15} (see also the final remark in \cite{p19}).
    \item All functors are implicitely derived so that, for example, given a proper morphism of schemes of finite Tor-dimension $f:X\rightarrow Y$ we will write $f_*$ instead of $\textup{R}f_*$. Similarly, all fiber products of schemes will be taken in the $\oo$-category of derived schemes.
    \item Whenever we are given a morphism of (derived) schemes $f:X\rightarrow Y$ and an object $\Ec_Y \in \qcoh{Y}$ of the ($\oo$-enhancement of the) derived category of $Y$, we will write $\Ec_X$ instead of $f^*\Ec_Y$.
    \item We will work with dg categories up to Morita equivalence. The $\oo$-category of ($S$-linear) dg categories up to Morita equivalences will be denoted by $\dgcatm{S}$.
\end{itemize}
\begin{center}
    \textit{Acknowledgements}
\end{center}
I would like to thank E.~Segal, B.~Toën and G.~Vezzosi for many useful conversations and remarks on the subject of this paper.

This project has received funding from the European Research Council (ERC) under the European Union Horizon 2020 research and innovation programme (grant agreement No. 725010).

\section{Cohomological operators}
\begin{context}\label{contesto di base}
Let $A$ be a local regular noetherian ring of finite Krull dimension, whose prime spectrum will be denoted by $S$. Let $Y$ be a flat $S$ scheme of finite type. Moreover, let $\Ec_Y$ be a vector bundle of rank $r$ on $Y$ and $\Lc_Y$ be a line bundle on $Y$. We will denote by $\Ec_Y^+$ the vector bundle $\Ec_Y\oplus \Lc_Y$ on $Y$.
\end{context}

Recall the definition of twisted LG models (see \cite{p20}):
\begin{defn}
Let $\tlgmpm{Y}$ denote the (ordinary) symmetric monoidal category of \ital{twisted LG models over} $(Y,\Ec_Y^+)$. It is defined as follows:
\begin{itemize}
    \item objects are pairs $(X,p^+)$, where $f:X\rightarrow Y$ is a flat $Y$-scheme of finite type (the structure morphism is tacit in our notation) and $p^+$ is a global section of $\Ec^+_X$. We will use the notation $(p,s)$ to denote the global section corresponding to $p^+$ under the isomorphism $\Hu^0(X,\Ec_X^+)\simeq \Hu^0(X,\Ec_X)\oplus \Hu^0(X,\Lc_X)$;
    
    \item given two objects $(X_{k},p^+_k)$, $k\in \{1,2\}$, a morphism $\phi:(X_1,p^+_1)\rightarrow (X_2,p^+_2)$ is a morphisms $\phi:X_1\rightarrow X_2$ of $Y$-schemes such that $\phi^*(p_2^+)=p_1^+$;
    
    \item identities and composition are clear;
    
    \item the symmetric monoidal structure, labelled 
          \begin{equation*}
              \boxplus : \tlgmp{Y}\times \tlgmp{Y}\rightarrow \tlgmp{Y},
          \end{equation*}
          is defined by the formula
          \begin{equation*}
              (X_1,p^+_1)\boxplus (X_2,p^+_2)= (X_1\times_Y X_2, p_1^+\boxplus p_2^+).
          \end{equation*}
          Here $p_1^+\boxplus p_2^+$ stands for the global section $pr_1^*(p_1^+)+pr_2^*(p_2^+)\in \Hu^0(X_1\times_Y X_2,\Ec^+_{X_1\times_Y X_2})$, where $pr_k:X_1\times_Y X_2\rightarrow X_k$ denotes the canonical projection, $k\in \{1,2\}$. It is straightforward that $\boxplus$ defines a symmetric monoidal structure, whose unit is the object $(Y,0)$.
\end{itemize}
\end{defn}
\begin{rmk}\label{functor Q}
Let $(Y,\Ec^+_Y)$ be as above and consider the symmetric monoidal category $\tlgmrm{Y}$ of twisted LG models over $(Y,\Ec_Y)$. There is an obvious symmetric monoidal functor  
   \begin{equation*}
      Q^{\boxplus}: \tlgmpm{Y}\rightarrow \tlgmrm{Y}
   \end{equation*}
   \begin{equation*}
       (X,p^+)\mapsto (X,p).
   \end{equation*}
Notice that $Q^{\boxplus}$ admits a section
   \begin{equation*}
       \iota : \tlgmrm{Y} \rightarrow \tlgmpm{Y}
   \end{equation*}
   \begin{equation*}
       (X,p)\mapsto(X,(p,0)).
   \end{equation*}
\end{rmk}
\begin{notation}
To any twisted LG model $(X,p^+)$ over $(Y,\Ec^+_Y)$ we can associate a (derived) scheme, namely the (derived) zero locus of $p^+$. We will denote it by $\Vu(p^+)$. Similarly, we will denote by $\Vu(p)$ the (derived) zero locus of $p$. 
\end{notation}
Recall the following construction of \cite[\S 6]{p20}. Let $(Z,\Vc_Z)$ be a proper flat scheme over $S$ with a fixed vector bundle. Let $\Pbb(\Vc_Z^{\vee})$ denote the projective bundle $Proj_Z(Sym_{\Oc_Z}(\Vc_Z))$ and let $\Oc(1)$ denote the twisting sheaf on it.
There exists a symmetric monoidal functor
\begin{equation*}
    \Xi^{\boxplus}: \textup{LG}_{(Z,\Vc_Z)}^{\boxplus}\rightarrow \tlgmprojm{\Vc_Z^{\vee}}
\end{equation*}
defined on objects by the formula
\begin{equation*}
    \Xi(X,s)=(\Pbb(\Vc_X^{\vee}),W_{s}).
\end{equation*}
Recall that 
\begin{equation*}
    W_{s}:\Oc_{\Pbb(\Vc^{\vee}_X)}=\widetilde{Sym_{\Oc_X}(\Vc_X)}\rightarrow \widetilde{Sym_{\Oc_X}(\Vc_X)}(1)=\Oc(1)
\end{equation*}
is the global section of $\Oc(1)$ induced by the morphism $\Oc_X\rightarrow \Vc_X$.

There are lax monoidal $\oo$-functors
\begin{equation*}
    \sing{\bullet,\bullet}, \sing{\Pbb(\Vc_{\bullet}^{\vee}),W_{\bullet}}: \textup{LG}_{(Z,\Vc_Z)}^{\textup{op},\boxplus} \rightarrow \dgcatmo{S}
\end{equation*}
which at the level of objects are defined by the assignments
\begin{equation*}
    (X,s)\mapsto \sing{X,p^+}, \hspace{0.5cm} (X,s)\mapsto \sing{\Pbb(\Vc_X^{\vee}),W_{s}}
\end{equation*}
respectively.
\begin{notation}
For us, given a pair $(Z,s)$, where $s$ is a global section of some vector bundle over $Z$, $\sing{Z,s}$ denotes the dg category of relative singularities of $i:\Vu(s)\hookrightarrow Z$, i.e. the homotopy fiber (computed in $\dgcatm{S}$) of the dg functor
\begin{equation*}
    i_*:\sing{\Vu(s)}\rightarrow \sing{Z}.
\end{equation*}
This dg functor is well defined as $i$ is an lci closed embedding of derived schemes. In particular, it preserves both perfect complexes and complexes with coherent cohomology (\cite{gr17}, \cite{to12}).
It is well known (see \cite{brtv18}, \cite{p20}) that 
\begin{equation*}
    \sing{Z,s}\simeq \frac{\cohbp{\Vu(s)}{Z}}{\perf{\Vu(s)}},
\end{equation*}
where $\cohbp{\Vu(s)}{Z}$ denotes the full sub dg category of $\cohb{\Vu(s)}$ spanned by those complexes $\Mcc$ such that $i_*(\Mcc)\in \perf{Z}$.
\end{notation}
Notice that since both functors are lax monoidal, they factor through $\md_{\sing{Z,0}}(\dgcatmo{S})$ and $\md_{\sing{\Pbb(\Vc_Z^{\vee}),0}}(\dgcatmo{S})$ respectively.
These two assignements are compatible in the following way:
\begin{thm}[\cite{or06}, \cite{bw15}, \cite{p20}]\label{reduction of codimension for sing}
There exists a lax monoidal $\oo$-natural trasformation of $\oo$-functors
\begin{equation*}
    \Upsilon^{\otimes}:\sing{\bullet,\bullet}\rightarrow \sing{\Pbb(\Vc_{\bullet}^{\vee}),W_{\bullet}} :\textup{LG}_{(Z,\Vc_Z)}^{\textup{op},\boxplus} \rightarrow \md_{\sing{Z,0}}(\dgcatmo{S})
\end{equation*}
inducing an equivalence
\begin{equation*}
    \Upsilon_{(X,s)}:\sing{X,s}\xrightarrow{\sim} \sing{\Pbb(\Vc_X^{\vee}),W_s}
\end{equation*}
for any LG model $(X,s)$ where $s$ is a regular global section.
\end{thm}

\begin{rmk}
We strongly believe that the regularity assumption imposed on $s$ is superfluous, at least if one is willing to consider the derived zero locus of $s$ instead of the classical one. 
\end{rmk}

In particular, this applies to $(Z,\Vc_Z)=(Y,\Ec_Y^+)$ and $(Z,\Vc_Z)=(Y,\Ec_Y)$.

With the same notation as above, let $\Vu(W_{p^+})$ and $\Vu(W_p)$ denote the zero loci of the global sections $W_{p^+}$ and $W_p$ respectively. They fit in the following commutative diagram of (possibly derived) $S$ schemes:
\begin{equation}\label{main diagram}
 \begindc{\commdiag}[20]
    \obj(-60,30)[1]{$\Pbb(\Ec^{+,\vee}_{\Vu(p^+)})$}
    \obj(0,30)[2]{$\Vu(W_{p^+})$}
    \obj(60,30)[3]{$\Pbb(\Ec_X^{+,\vee})$}
    \obj(-60,10)[4]{$\Vu(p^+)$}
    \obj(30,0)[5]{$X$}
    \obj(-60,-10)[6]{$\Vu(p)$}
    \obj(-60,-30)[7]{$\Pbb(\Ec_{\Vu(p)}^{\vee})$}
    \obj(0,-30)[8]{$\Vu(W_p)$}
    \obj(60,-30)[9]{$\Pbb(\Ec_X^{\vee})$}
    \mor{1}{2}{$k$}
    \mor{2}{3}{$$}
    \mor{1}{4}{$q$}[\atright,\solidarrow]
    \mor{4}{6}{$i$}[\atright,\solidarrow]
    \mor(-55,10)(3,10){$$}[\atright, \solidline]
    \cmor((2,10)(15,10)(28,10)(30,8)(30,4)) \pdown(-30,10){$$}
    \mor(-55,-10)(3,-10){$$}[\atright, \solidline]
    \cmor((2,-10)(15,-10)(28,-10)(30,-8)(30,-4)) \pup(-30,-10){$$}
    \mor{7}{6}{$p$}
    \mor{7}{8}{$j$}
    \mor{8}{9}{$$}
    \mor{8}{2}{$r$}
    \mor{9}{3}{$$}
    \mor{3}{5}{$$}
    \mor{9}{5}{$$}
 \enddc
\end{equation}

In the diagram above, $\Pbb(\Ec_X^{\vee})\rightarrow \Pbb(\Ec_X^{+,\vee})$ is the closed embedding induced by the surjection $\Ec_X^{+}\simeq \Ec_X\oplus \Lc_X\rightarrow \Ec_X$. Notice that it is an lci closed embedding, as locally it is of the form $\proj{r-1}{U}\rightarrow \proj{r}{U}$. Also notice that $\Vu(W_{p})\simeq \Vu(W_{p^+})\times_{\Pbb(\Ec_X^{+,\vee})}\Pbb(\Ec_X^{\vee})$. In particular, $r:\Vu(W_p)\rightarrow \Vu(W_{p^+})$ is an lci morphism as well and therefore it induces a dg functor
\begin{equation*}
    r_*:\sing{\Pbb(\Ec_X^{\vee}),W_p}\rightarrow \sing{\Pbb(\Ec_X^{+,\vee}),W_{p^+}}
\end{equation*}
at the level of dg categories of relative singularities. 
\begin{lmm}
The dg functor $r_*:\sing{\Vu(p),W_p}\rightarrow \sing{\Vu(p^+),W_{p^+}}$ admits a right adjoint $r^!: \sing{\Vu(p^+),W_{p^+}} \rightarrow \sing{\Vu(p),W_p}$ induced by the !-pullback.
\end{lmm}
\begin{proof}
This is a standard verification using the properties of the !-pullback. Consider the pullback square
\begin{equation*}
    \begindc{\commdiag}[20]
    \obj(-18,15)[1]{$\Vu(W_{p})$}
    \obj(18,15)[2]{$\Pbb(\Ec_X^{\vee})$}
    \obj(-18,-15)[3]{$\Vu(W_{p^+})$}
    \obj(18,-15)[4]{$\Pbb(\Ec_X^{+,\vee}).$}
    \mor{1}{2}{$\alpha$}
    \mor{1}{3}{$r$}[\atright,\solidarrow]
    \mor{2}{4}{$\tilde{r}$}
    \mor{3}{4}{$\tilde{\alpha}$}[\atright,\solidarrow]
    \enddc
\end{equation*}
Since $r$ is lci and, in particular, of finite Tor-dimension, we have a well defined dg functor $r^*:\cohb{\Vu(W_{p^+})}\rightarrow \cohb{\Vu(p)}$. Its restriction to $\cohbp{\Vu(W_{p^+})}{\Pbb(\Ec^{+,\vee}_X)}$ factors through $\cohbp{\Vu(W_p)}{\Pbb(\Ec_X^{\vee})}\subseteq \cohb{\Vu(W_p)}$ as $r^*\alpha_*\simeq \tilde{r}^*\tilde{\alpha}_*$. Recall that we have an equivalence of functors $r^!(-)\simeq r^*(-)\otimes \omega_{\Vu(W_p)/\Vu(W_{p^+})}$ (see, e.g., \cite[Corollary 6.4.2.7]{sag}), where $\omega_{\Vu(W_p)/\Vu(W_{p^+})}$ denotes the relative dualizing sheaf. Notice that $\omega_{\Vu(W_p)/\Vu(W_{p^+})}\simeq \alpha^*\omega_{\Pbb(\Ec_X^{\vee})/\Pbb(\Ec_X^{+,\vee})}$ is a perfect complex. In particular, $r^!$ is well defined at the level of complexes with coherent bounded cohomology sheaves. Then we have a well defined dg functor
\begin{equation*}
    r^!:\cohbp{\Vu(W_{p+})}{\Pbb(\Ec^{+,\vee}_X)}\rightarrow \cohbp{\Vu(W_{p})}{\Pbb(\Ec_X^{\vee})}
\end{equation*}
as $\alpha_*(r^*(\Fc)\otimes \omega_{\Vu(W_p)/\Vu(W_{p^+})})\simeq \alpha_*r^*(\Fc)\otimes \omega_{\Pbb(\Ec_X^{\vee})/\Pbb(\Ec^{+,\vee}_X)}$. Since it is clear that $r^!$ preserves perfect complexes, we have a well defined dg functor
\begin{equation*}
    r^!: \sing{\Vu(p^+),W_{p^+}} \rightarrow \sing{\Vu(p),W_p}.
\end{equation*}
It is clear that it is the right adjoint to $r_*:\sing{\Vu(p),W_p} \rightarrow \sing{\Vu(p^+),W_{p^+}}$ (before taking quotients, $r^!$ is right adjoint to $r_*$).
\end{proof}
Let us now focus on the closed embedding $i:\Vu(p^+)\rightarrow \Vu(p)$. This is also an lci closed embedding of (derived) schemes and therefore the \textasteriskcentered-pullback and the \textasteriskcentered-pushforward induce an adjunction at the level of the dg categories of singularities (see \cite{ep15}, \cite{brtv18}, \cite{p20}). We will denote this adjunction by
\begin{equation*}
    i^*: \sing{X,p} \rightleftarrows \sing{X,p^+} :i_* .
\end{equation*}
It is natural to ask whether the adjunction $(i^*,i_*)$ identifies with $(r_*,r^!)$ under the equivalences
\begin{equation*}
    \sing{X,p^+}\xrightarrow{\sim} \sing{\Pbb(\Ec^{+,\vee}_X),W_{p^+}}, \hspace{0.5cm} \sing{X,p}\xrightarrow{\sim} \sing{\Pbb(\Ec_X^{\vee}),W_{p}}.
\end{equation*}
We will now show that this is indeed the case up to a twisted shift.
\begin{rmk}
This question is closely related with the comparison between Eisenbud operators and operators induced by $T_i\in \Hu^0(\proj{r-1}{B},\Oc(1))$ addressed by the authors in \cite{bw15}. We will come back to this point later.
\end{rmk}
\begin{thm}\label{equivalence adjunction}
Let $(X,p^+)$ be an LG model over $(Y,\Ec_Y^+)$ such that $p^+$ is regular. There is an equivalence of adjunctions 
\begin{equation*}
    \begindc{\commdiag}[20]
       \obj(-70,15)[1]{$i^*\bigl ((-)\otimes \Lc_{\Vu(p)}  \bigr )[-1]:\sing{X,p}$}
       \obj(70,15)[2]{$\sing{X,p^+}:\bigl (i_*(-)\otimes \Lc_{\Vu(p)}^{\vee}\bigr )[1]$}
       \obj(-70,-15)[3]{$r_*:\sing{\Pbb(\Ec_X),W_{p}}$}
       \obj(70,-15)[4]{$\sing{\Pbb(\Ec^+_X),W_{p^+}}:r^!$}
       \mor(-30,18)(30,18){$$}
       \mor(30,12)(-30,12){$$}
       \mor(30,-18)(-30,-18){$$}
       \mor(-30,-12)(30,-12){$$}
       \mor{1}{3}{$\sim$}[\atright,\solidarrow]
       \mor{2}{4}{$\sim$}
    \enddc
\end{equation*}
where the vertical dg functors are $\Upsilon_{(X,p)}$ and $\Upsilon_{(X,p^+)}$ respectively, i.e.
\begin{equation*}
    r_*\circ \Upsilon_{(X,p)}\simeq \Upsilon_{(X,p^+)}\circ i^*\bigl ((-)\otimes \Lc_{\Vu(p)}  \bigr )[-1], \hspace{0.5cm} r^!\circ \Upsilon_{(X,p^+)}\simeq \Upsilon_{(X,p)}\circ \bigl (i_*(-)\otimes \Lc_{\Vu(p)}^{\vee}\bigr )[1].
\end{equation*}
\end{thm}
\begin{proof}
It is clear that it suffices to show one of the equivalences in the statement of the theorem. For example, if we assume that $ r_*\circ \Upsilon_{(X,p)}\simeq \Upsilon_{(X,p^+)}\circ i^*\bigl ((-)\otimes \Lc_{\Vu(p)}  \bigr )[-1]$, we get the other equivalence by noticing that, by uniqueness (up to equivalence) of right adjoints, $\Upsilon_{(X,p)}^{-1}\circ r^!\simeq \bigl (i_*(-)\otimes \Lc_{\Vu(p)}^{\vee}\bigr )[1]\circ \Upsilon_{(X,p^+)}^{-1}$, i.e. $r^!\circ \Upsilon_{(X,p^+)}\simeq \Upsilon_{(X,p)}\circ \bigl (i_*(-)\otimes \Lc_{\Vu(p)}^{\vee}\bigr )[1]$. The other implication is analogous. 

We will show that $\Upsilon_{(X,p)}^{-1}\circ r^!\simeq \bigl (i_*(-)\otimes \Lc_{\Vu(p)}^{\vee}\bigr )[1] \circ \Upsilon_{(X,p^+)}^{-1}$. Recall that $\Upsilon_{(X,p)}\simeq j_*p^*$ and that $\Upsilon_{(X,p^+)}\simeq k_*q^*$. In particular, $\Upsilon_{(X,p)}^{-1}\simeq p_*j^!$ and $\Upsilon_{(X,p^+)}^{-1}\simeq q_*k^!$. We are therefore claiming that
\begin{equation*}
    p_*j^!r^!\simeq \bigl (i_*(-)\otimes \Lc_{\Vu(p)}^{\vee}\bigr )q_*k^![1].
\end{equation*}
Let us compose both functors with $\Upsilon_{(X,p^+)}=k_*q^*$ on the right. On the right hand side we immediately see that we obtain a functor equivalent to $(i_*(-)\otimes \Lc_{\Vu(p)}^{\vee}\bigr )[1]$. 

Consider the following diagram
\begin{equation*}
    \begindc{\commdiag}[20]
    \obj(-50,20)[1]{$\Vu(p^+)$}
    \obj(0,20)[2]{$\Pbb(\Ec^{+,\vee}_{\Vu(p^+)})$}
    \obj(50,20)[3]{$\Vu(W_{p^+})$}
    \obj(0,0)[4]{$\Pbb(\Ec^{\vee}_{\Vu(p^+)})$}
    \obj(-50,-20)[5]{$\Vu(p)$}
    \obj(0,-20)[6]{$\Pbb(\Ec^{\vee}_{\Vu(p)})$}
    \obj(50,-20)[7]{$\Vu(W_p)$}
    \mor{2}{1}{$q$}
    \mor{2}{3}{$k$}
    \mor{4}{1}{$q\circ \alpha$}
    \mor{4}{2}{$\alpha$}
    \mor{1}{5}{$i$}
    \mor{4}{6}{$\beta$}
    \mor{6}{5}{$p$}
    \mor{4}{7}{$\beta \circ j$}
    \mor{6}{7}{$j$}[\atright,\solidarrow]
    \mor{7}{3}{$r$}[\atright,\solidarrow]
    \enddc
\end{equation*}
and notice that $\Vu(p^+)\times_{\Vu(p)}\Pbb(\Ec^{\vee}_{\Vu(p)})\simeq \Pbb(\Ec^{\vee}_{\Vu(p^+)})\simeq \Pbb(\Ec^{+,\vee}_{\Vu(p^+)})\times_{\Vu(W_{p^+})}\Vu(W_p)$. Then we have the following chain of equivalences
\begin{equation*}
    p_*j^!r^!k_*q^*\underbracket{\simeq}_{r^!k_*\simeq j_*\beta_*\alpha^!} p_*j^!j_*\beta_*\alpha^!q^*\underbracket{\simeq}_{\alpha^!\simeq \alpha^*(-)\otimes  \omega_{\Pbb(\Ec^{\vee}_{\Vu(p^+)})/\Pbb(\Ec^{+,\vee}_{\Vu(p^+)})}} p_*j^!j_*\beta_*(\alpha^*q^*(-)\otimes \omega_{\Pbb(\Ec^{\vee}_{\Vu(p^+)})/\Pbb(\Ec^{+,\vee}_{\Vu(p^+)})}).
\end{equation*}
The ideal sheaf $\Is \subseteq \Oc_{\Pbb(\Ec_X^{+,\vee})}$ defining $u:\Pbb(\Ec_X^{\vee})\rightarrow \Pbb(\Ec_X^{+,\vee})$ is 
\begin{equation*}
    \Oc_{\Pbb(\Ec_X^{+,\vee})}(-1)\otimes_{\Oc_X}\Lc_X\simeq  \Oc_{\Pbb(\Ec_X^{+,\vee})}(-1)\otimes_{\Oc_{\Pbb(\Ec_X^{+,\vee})}}\Lc_{\Pbb(\Ec_X^{+,\vee})}, 
\end{equation*} 
from which we see that
\begin{equation*}
    \frac{\Is}{\Is^2}\simeq \Oc_{\Pbb(\Ec_X^{\vee})}(-1)\otimes_{\Oc_{\Pbb(\Ec_X^{\vee})}}\Lc_{\Pbb(\Ec_X^{\vee})}.
\end{equation*}
Then,
\begin{equation*}
    \omega_{\Pbb(\Ec_X^{\vee})/\Pbb(\Ec_X^{+,\vee})}\simeq \Ns_{\Pbb(\Ec_X^{\vee})/\Pbb(\Ec_X^{+,\vee})}[-1]=\Bigl ( \frac{\Is}{\Is^2} \Bigr )^{\vee}[-1]\simeq (\Oc_{\Pbb(\Ec_X^{\vee})}(1)\otimes \Lc_{\Pbb(\Ec_X^{\vee})}^{\vee})[-1].
\end{equation*}
As the squares in the diagram
\begin{equation*}
    \begindc{\commdiag}[20]
    \obj(-40,10)[1]{$\Pbb(\Ec_{\Vu(p^+)}^{+,\vee})$}
    \obj(0,10)[2]{$\Vu(W_{p^+})$}
    \obj(40,10)[3]{$\Pbb(\Ec_{X}^{+,\vee})$}
    \obj(-40,-10)[4]{$\Pbb(\Ec_{\Vu(p^+)}^{\vee})$}
    \obj(0,-10)[5]{$\Vu(W_p)$}
    \obj(40,-10)[6]{$\Pbb(\Ec_{X}^{\vee})$}
    \mor{1}{2}{$k$}
    \mor{2}{3}{$$}
    \mor{4}{1}{$\alpha$}
    \mor{5}{2}{$r$}
    \mor{6}{3}{$$}
    \mor{4}{5}{$\beta \circ j$}
    \mor{5}{6}{$u$}
    \enddc
\end{equation*}

are cartesian, the equivalences 
\begin{equation*}
    \omega_{\Pbb(\Ec_{\Vu(p^+)}^{\vee})/\Pbb(\Ec^{+,\vee}_{\Vu(p^+)})}\simeq (\Oc_{\Pbb(\Ec_{\Vu(p^+)}^{\vee})}(1)\otimes \Lc_{\Pbb(\Ec^{\vee}_{\Vu(p^+)})}^{\vee})[-1]\simeq \beta^*j^*(\Oc_{\Vu(W_p)}(1)\otimes \Lc_{\Vu(W_p)}^{\vee})[-1]
\end{equation*} 
follow from \cite[Remark 6.4.2.6]{sag}.
Therefore, we can continue the chain of equivalences
\begin{equation*}
    p_*j^!j_*\beta_*(\alpha^*q^*(-)\otimes \omega_{\Pbb(\Ec^{\vee}_{\Vu(p^+)})/\Pbb(\Ec^{+,\vee}_{\Vu(p^{+})})})\simeq p_*j^!j_*\beta_*\bigl (\alpha^*q^*(-)\otimes \beta^*j^*(\Oc_{\Vu(W_p)}(1)\otimes_{\Oc_{\Vu(W_p)}}\Lc_{\Vu(W_p)}^{\vee})[-1]\bigr )
\end{equation*}    
\begin{equation*}
   \simeq p_*j^!\bigl (j_*\beta_*\alpha^*q^*(-)\otimes \Oc_{\Vu(W_p)}(1)\otimes \Lc_{\Vu(W_p)}^{\vee}\bigr )[-1] \simeq p_*j^!\bigl (j_*\beta_*\alpha^*q^*(-)\otimes \Lc_{\Vu(W_p)}^{\vee}\bigr )[1],
\end{equation*}
where for the last equivalences we have used the projection formula and the fact that
\begin{equation*}
    (-)\otimes \Oc_{\Vu(W_p)}(1)[-1] \simeq (-)[1]
\end{equation*}
in $\sing{\Pbb(\Ec_X^{\vee}),W_p}$ (i.e., $\sing{\Pbb(\Ec_X^{\vee}),W_p}$ is twisted $2$-periodic, see e.g. \cite[Theorem 2.7]{bw15}).
Finally, notice that
\begin{equation*}
    p_*j^!\bigl (j_*\beta_*\alpha^*q^*(-)\otimes \Lc_{\Vu(W_p)}^{\vee}\bigr )[1]\simeq p_*j^!j_*\bigl (\beta_*\alpha^*q^*(-)\otimes \Lc_{\Pbb(\Ec^{\vee}_{\Vu(p)})}^{\vee}\bigr )[1]
\end{equation*}
\begin{equation*}
    \underbracket{\simeq}_{\beta_*\alpha^*q^*\simeq p^*i_*}p_*j^!j_*\bigl (p^*i_*(-)\otimes \Lc_{\Pbb(\Ec^{\vee}_{\Vu(p)})}^{\vee}\bigr )[1]\simeq p_*j^!j_*p^*\bigl (i_*(-)\otimes \Lc_{\Vu(p)}^{\vee}\bigr )[1]
\end{equation*}
\begin{equation*}
    \underbracket{\simeq}_{p_*j^!j_*p^*\simeq id_{\sing{X,p}}} \bigl (i_*(-)\otimes \Lc_{\Vu(p)}^{\vee}\bigr )[1].
\end{equation*}
Since $k_*q^*$ is an equivalence, the proof is complete.
\end{proof}
It has been known since when they were introduced by Eisenbud in \cite{ei80} that matrix factorizations are naturally endowed with certain \textit{cohomological operators}, known nowadays as \textit{Eisenbud operators}. They were defined in the following setup: consider a local regular noetherian ring of finite Krull dimension $Q$ and let $f=(f_1,\dots,f_n)\in Q^n$ be a regular sequence on it. Set $R=Q/(f)$. Then Eisenbud defined in \textit{loc. cit.} $n$ natural trasformations in $\textup{Fun}(\sing{R},\sing{R})$
\begin{equation*}
    \varepsilon_j:id_{\sing{R}}\rightarrow (-)[2], \hspace{0.5cm} j\in \{1,\dots,n\}.
\end{equation*}
\begin{rmk}
In the special case $n=1$, $\varepsilon_1=\varepsilon$ is an equivalence and is indeed the one giving $\sing{R}\simeq \mf{Q}{f}$ its natural $2$-periodic structure.
\end{rmk}
In \cite{bw15}, the authors provide a geometric description of Eisenbud operators, using the equivalence provided by Theorem \ref{reduction of codimension for sing}. Notice that they don't require $Q$ to be regular.
They show (see \cite[Theorem 3.2]{bw15}) that under the equivalence
\begin{equation*}
    \sing{Q,f}\simeq \sing{\proj{n-1}{Q},W},
\end{equation*}
where $W=\sum_{j=1}^nf_j\cdot T_j$ ($T_j$'s being the local coordinates of $\proj{n-1}{Q}$), Eisenbud operators correspond to the natural trasformations
\begin{equation*}
    id_{\sing{\proj{n-1}{Q},W}}\xrightarrow{T_j} (-)\otimes \Oc(1)\simeq (-)[2]
\end{equation*}
induced by multiplication along the local coordinates of the projective space.

One of the main purposes of this note is to show that it might convenient to consider Eisenbud operators in more general contexts then the one quickly reminded above. For this we introduce the following
\begin{defn}\label{defn E-operator}
Let $(Y,\Ec_Y^+)$ be as in Context \ref{contesto di base} and let $(X,p^+)$ be a LG model over it. Using the same notation as in diagram \ref{main diagram}, define $\varepsilon$ to be the natural trasformation
\begin{equation*}
    id_{\sing{X,p^+}}\xrightarrow{\varepsilon} (-)\otimes \Lc_{\Vu(p^+)}^{\vee}[2]
\end{equation*}
induced by the counit $i^*i_*\rightarrow  id_{\sing{X,p^+}}$. In particular, for any object $\Bc$ in $\sing{X,p^+}$, by definition, we have a fiber-cofiber sequence
\begin{equation*}
    i^*i_*\Bc \rightarrow \Bc \xrightarrow{\varepsilon} \Bc \otimes \Lc_{\Vu(p^+)}^{\vee}[2].
\end{equation*}
\begin{comment}
\red{add lemma computing the cofiber of $i^*i_*\Bc \rightarrow \Bc$}.
\end{comment}
We will refer to $\varepsilon$ as the \textit{Eisenbud operator} on $\sing{X,p^+}$ (E-operator for short). 
\end{defn}

Similarly, we have an operator on the other side of the equivalence in Theorem \ref{reduction of codimension for sing}.

\begin{defn}\label{defn BW operator}
Let $(Y,\Ec_Y^+)$ be as in Context \ref{contesto di base} and let $(X,p^+)$ be a LG model over it. Using the same notation as in diagram \ref{main diagram}, define $\chi$ to be the natural trasformation
\begin{equation*}
    id_{\sing{\Pbb(\Ec_X^+),W_{p^+}}}\xrightarrow{\chi} (-)\otimes \Oc_{\Vu(W_{p^+})}(1) \otimes \Lc_{\Vu(W_{p^+})}^{\vee}
\end{equation*}
induced by the counit $r_*r^!\rightarrow id_{\sing{\Pbb(\Ec_X^+),W_{p^+}}}$. In particular, for any object $\Cc$ in $\sing{\Pbb(\Ec_X^{+,\vee}),W_{p^+}}$, we have a fiber-cofiber sequence
\begin{equation*}
    r_*r^!\Cc \rightarrow \Cc \xrightarrow{\chi} \Cc\otimes \Oc_{\Vu(W_{p^+})}(1) \otimes \Lc_{\Vu(W_{p^+})}^{\vee}.
\end{equation*}
\begin{comment}
\red{add lemma computing the cofiber of $ r_*r^!\Cc \rightarrow \Cc$ }.
\end{comment}
We will refer to $\chi$ as the \textit{Burke-Walker operator} on $\sing{\Pbb(\Ec^{+,\vee}_X),W_{p^+}}$ (BW-operator for short).
\end{defn}
As an immediate consequence of Theorem \ref{equivalence adjunction}, we have
\begin{cor}\label{E-operator = BW-operator}
    Let $(X,p^+)$ be an LG model over $(Y,\Ec_Y^+)$ such that $p^+$ is regular. Then the E-operator corresponds to the BW-operator under the equivalence $\Upsilon_{(X,p^+)}$.
\end{cor}
\begin{ex}
Let $Q,f$ and $R$ be as in the discussion preceding Definition \ref{defn E-operator}. Let $Y=X=S=Spec(Q)$, $\Ec_Y=B^{n-1}$, $\Lc_Y=B$, $p=(f_1,\dots,f_{n-1})$ and $s=f_n$ (so that $p^+=f\in B^n=\Ec_X^+$). Then $\Vu(p^+)=Spec(R)$ and $\Vu(p)=Spec(Q/(f_1,\dots,f_{n-1}))$. The morphism $\Pbb(\Ec^{\vee}_X)\rightarrow \Pbb(\Ec_X^{+,\vee})$ identifies with the hyperplane $\proj{n-2}{Q}\hookrightarrow \proj{n-1}{Q}$ cut by the global section $T_n\in \Hu^0(\proj{n-1}{Q},\Oc(1))$. In particular, $\Vu(W_p)=\Vu(W_{p^+},T_n) \hookrightarrow  \Vu(W_{p^+})$. In this case, for any $\Cc \in \sing{\proj{n-1}{Q},W_{p^+}}$, the cofiber of the morphism $r_*r^!\Cc \rightarrow \Cc$ identifies with $\Cc\xrightarrow{T_n} \Cc \otimes \Oc(1)$ and therefore the BW-operator $\chi$ agrees with the natural trasformation defined on $\mf{\proj{n-1}{Q},\Oc(1)}{W_{p^+}}\simeq \sing{\proj{n-1}{Q},W_{p^+}}$ in \cite{bw15}. Moreover, as a consequence of \cite[Theorem 3.2]{bw15} and of Corollary \ref{E-operator = BW-operator}, we obtain that the E-operator $\varepsilon$ agrees with the operator defined in \cite{ei80}. This fact should reassure the reader that our terminology is consistent with the existing literature on the subject.
\end{ex}
\begin{rmk}
Consider the same notation as in the previous example. It is known since their introduction in \cite{ei80} that the Eisenbud operators $\varepsilon_i$ commute with each other. This can be most easily seen using Burke-Walker's theorem \cite[Theorem 3.2]{bw15}. After our identification of Eisenbud operators with the cofiber of the counit morphism 
\begin{equation*}
    i_{0j}^*i_{0j*}\Bc \rightarrow \Bc,
\end{equation*}
where $i_{0j}:Spec(Q/f)\hookrightarrow Spec(Q/(f_1,\dots,\hat{f_j},\dots,f_n))$, this can also be seen by  looking at the following commutative diagram
\begin{equation*}
    \begindc{\commdiag}[20]
       \obj(-40,20)[1]{$i_{0k}^*i_{0k*}i_{0j}^*i_{0j*}F$}
       \obj(0,20)[2]{$i_{0j}^*i_{0j*}F$}
       \obj(40,20)[3]{$i_{0j}^*i_{0j*}F[2]$}
       \obj(-40,0)[4]{$i_{0k}^*i_{0k*}F$}
       \obj(0,0)[5]{$F$}
       \obj(40,0)[6]{$F[2]$}
       \obj(-40,-20)[7]{$i_{0k}^*i_{0k*}F[2]$}
       \obj(0,-20)[8]{$F[2]$}
       \obj(40,-20)[9]{$F[4].$}
       \mor{1}{2}{$$}
       \mor{2}{3}{$\varepsilon_{k}$}
       \mor{1}{4}{$$}
       \mor{2}{5}{$$}
       \mor{3}{6}{$$}
       \mor{4}{5}{$$}
       \mor{5}{6}{$\varepsilon_{k}$}
       \mor{4}{7}{$\varepsilon_{j}$}
       \mor{5}{8}{$\varepsilon_{j}$}
       \mor{6}{9}{$\varepsilon_{j}$}
       \mor{7}{8}{$$}
       \mor{8}{9}{$\varepsilon_{k}$}
    \enddc
\end{equation*}
\end{rmk}
\section{dg category of coherent relative singularities}
Let $(X,p^+)$ be an LG model over $(Y,\Ec_Y^+)$. In this section we will investigate torsion objects for the E-operator (resp. BW-operator) on $\sing{X,p^+}$ (resp. $\sing{\Pbb(\Ec_X^{+,\vee}),W_{p^+}}$). When $X$ and $p^+$ are regular, this will allow us to identify the dg category of coherent singularities of $i:\Vu(p^+)\hookrightarrow \Vu(p)$ with a certain dg category of absolute singularities.

Recall that there exists a lax monoidal $\oo$-functor
\begin{equation*}
    \singpfm{Y}: \tlgmpmop{Y}\rightarrow \dgcatmo{S}
\end{equation*}
which at the level of objects is defined by
\begin{equation*}
    (X,p^+)\mapsto \singpf{Y}(X,p^+)= Ker\bigl ( f_*:\sing{\Vu(p^+)}\rightarrow \sing{X} \bigr ) \simeq \frac{\cohbp{\Vu(p^+)}{X}}{\perf{\Vu(p^+)}},
\end{equation*}
where $f$ denotes the lci closed embedding $\Vu(p^+)\hookrightarrow X$.
Similarly, there exists a lax monoidal $\oo$-functor
\begin{equation*}
    \singfm{Y}: \tlgmrmop{Y}\rightarrow \dgcatmo{S}
\end{equation*}
which at the level of objects is defined by
\begin{equation*}
    (X,p)\mapsto \singf{Y}(X,p)= Ker\bigl ( g_*:\sing{\Vu(p)}\rightarrow \sing{X} \bigr ) \simeq \frac{\cohbp{\Vu(p)}{X}}{\perf{\Vu(p)}},
\end{equation*}
where $g$ is the lci closed embedding $\Vu(p)\hookrightarrow{X}$.
By composing $\singfm{Y}$ with $Q^{\boxplus}$ (see Remark \ref{functor Q}), we obtain a lax monoidal $\oo$-functor
\begin{equation*}
    \singfm{Y}\circ Q^{\boxplus}, \hspace{0.5cm} (X,p^+)\mapsto \singf{Y}(X,p).
\end{equation*}
\begin{rmk}
We apologize for the heavy notations $\singfm{Y}$ and $\singpfm{Y}$, but we think it would lead to confusion to denote with $\textup{\textbf{Sing}}^{\otimes}$ both the $\oo$-functor $\tlgmrmop{Y}\rightarrow \dgcatmo{S}$ and $\tlgmpmop{Y}\rightarrow \dgcatmo{S}$.
\end{rmk}
\begin{prop}\label{functoriality (i^*,i_*)}
Let $(X,p^+)$ be an LG model over $(Y,\Ec_Y^+)$ and let $i:\Vu(p^+)\rightarrow \Vu(p)$ denote the canonical lci closed embedding. The adjunction
\begin{equation*}
    i^*:\sing{X,p}\rightleftarrows \sing{X,p^+}:i_*
\end{equation*}
is functorial in $(X,p^+)$, i.e. we have two natural transformations
\begin{equation*}
\singf{Y}\circ Q\rightarrow \singpf{Y}, \hspace{0.5cm} \singpf{Y}\circ Q \rightarrow \singf{Y}    
\end{equation*}
that, for a fixed LG model $(X,p^+)$ over $(Y,\Ec_Y^+)$, determine the above mentioned adjunction.

In particular, the E-operator is functorial in $(X,p^+)$.
\end{prop}
\begin{proof}
By uniqueness (up to equivalence) of right adjoints, it suffices to show that the dg functor $i^*:\sing{X,p}\rightarrow \sing{X,p^+}$ is functorial in $(X,p^+)$. Recall that the lax monoidal $\oo$-functors $\singpfm{Y}$ and $\singfm{Y}$ are defined using strict models and Kan extensions. We refer the reader to \cite[\S3 and \S6]{p20} for details. In particular, we can reduce to consider the case where $Y$ and $X$ are affine. Let $(Y,\Ec_Y^+)$ correspond to $(A,E_A^+)$, $E_A^+\simeq E_A\oplus L_A$, where $E_A$ is a projective $A$ module of rank $r$ and $L_A$ is a line bundle. Let $(B,p^+)$ be an affine LG model over $(A,E_A^+)$. The dg category $\cohs{B,(E_B^+)^{\vee},p^+}$ is a strict model for the dg category $\cohbp{V(p^+)}{Spec(B)}$ (see \cite[Construction 6.1.6]{p20}), where $V(p^+)=Spec(K^s(B,(E_B^+)^{\vee},p^+))$ is the spectrum of the (simplicial version of the) Koszul algebra associated to $((E_B^+)^{\vee},p^+)$. Similarly, $\cohs{B,E_B^{\vee},p}$ is a strict model for the dg category $\cohbp{V(p)}{Spec(B)}$. Given a morphism $\phi:(B,p^+)\rightarrow (C,q^+)$ of affine LG models over $(A,E_A^+)$, i.e. a morphism $\phi:B\rightarrow C$ of $A$ algebras such that $p^+=(p,s)\mapsto q^+=(q,t)$, the pseudo-functoriality on $\textup{Coh}^s$ is induced by $-\otimes_BC$:
\begin{equation*}
    \cohs{B,(E_B^+)^{\vee},p^+}\rightarrow \cohs{C,(E_C^+)^{\vee},q^+}
\end{equation*}
\begin{equation*}
    F\mapsto F\otimes_B C.
\end{equation*}
The pseudo functoriality $\cohs{B,E_B^{\vee},p}\rightarrow \cohs{C,E_C^{\vee},q}$ is analogous. The pullback $i^*:\cohbp{\Vu(p)}{Spec(B)}\rightarrow \cohbp{\Vu(p^+)}{Spec(B)}$ is modelled by the following pseudo-functor:
\begin{equation*}
    \cohs{B,E_B^{\vee},p}\rightarrow \cohs{B,(E_B^+)^{\vee},p^+}
\end{equation*}
\begin{equation*}
    F\mapsto F\otimes_B K(L_B^{\vee},s).
\end{equation*}
Therefore, we need to show that the square
\begin{equation*}
    \begindc{\commdiag}[20]
    \obj(-40,20)[1]{$\cohs{B,E_B^{\vee},p}$}
    \obj(40,20)[2]{$\cohs{C,E_C^{\vee},q}$}
    \obj(-40,-20)[3]{$\cohs{B,(E_B^+)^{\vee},p^+}$}
    \obj(40,-20)[4]{$\cohs{C,(E_C^+)^{\vee},q^+}$}
    \mor{1}{2}{$-\otimes_BC$}
    \mor{1}{3}{$-\otimes_BK(L_B^{\vee},s)$}[\atright,\solidarrow]
    \mor{2}{4}{$-\otimes_CK(L_C^{\vee},t)$}
    \mor{3}{4}{$-\otimes_BC$}[\atright,\solidarrow]
    \enddc
\end{equation*}
commutes up to canonical equivalence. This follows immediately from the canonical equivalences
\begin{equation*}
    C\otimes_C-\simeq id,\hspace{0.5cm} K(L_C^{\vee},t)\simeq C\otimes_B K(L_B^{\vee},s).
\end{equation*}

The fact that functoriality of $i^*:\cohbp{\Vu(p)}{X}\rightarrow \cohbp{\Vu(p^+)}{X}$ in $(X,p^+)$ implies that of $i^*:\sing{X,p}\rightarrow \sing{X,p^+}$ is clear.

The last assertion of the statement is obvious.
\end{proof}
\begin{rmk}
Notice that the proof of the previous proposition shows the stronger assertion that the adjunction
\begin{equation*}
    i^*:\cohbp{\Vu(p)}{X}\rightleftarrows \cohbp{\Vu(p^+)}{X}:i_*
\end{equation*}
is functorial in $(X,p^+)$.
\end{rmk}
Our next aim will be to study $\varepsilon$-torsion objects inside $\sing{X,p^+}$.
\begin{defn}
Let $(X,p^+)$ be an LG model over $(Y,\Ec_Y^+)$. The dg category of $\varepsilon$-torsion objects in $\sing{X,p^+}$ is the full sub dg category 
\begin{equation*}
    \sing{X,p^+}^{\varepsilon-tors}:=\Bigl \{\Bc \in \sing{X,p^+}: \exists n\geq0 \text{ such that } \varepsilon^n\sim 0: \Bc \rightarrow \Bc \otimes \Lc_{\Vu(p^+)}^{\vee,\otimes n}[2n] \Bigr \}.
\end{equation*}
\end{defn}
\begin{prop}\label{E-torsion = image pullback}
Let $\ising{X,p}{i}$ denote the thick full sub dg category of $\sing{X,p^+}$ spanned by the image of the dg functor
\begin{equation*}
    i^*:\sing{X,p}\rightarrow \sing{X,p^+}.
\end{equation*}
Then,
\begin{equation*}
    \sing{X,p^+}^{\varepsilon-tors}\simeq \ising{X,p}{i}.
\end{equation*}
\end{prop}
\begin{proof}
It is clear that $\sing{X,p^+}^{\varepsilon-tors}\subseteq \sing{X,p^+}$ is a thick sub dg category. Therefore, it will suffice to show that the image of $i^*:\sing{X,p}\rightarrow \sing{X,p^+}$ is contained in $\sing{X,p^+}^{\varepsilon-tors}$. Let $\Pc$ be an object of $\sing{X,p}$. We will show that $i^*\Pc$ is a $\varepsilon$-torsion object. The counit morphism $i^*i_*i^*\Pc \xrightarrow{c}i^*\Pc$ admits a section $i^*\Pc \xrightarrow{u} i^*i_*i^*\Pc$, i.e. there exists a $2$-cell $\sigma$ providing an homotopy $id_{i^*\Pc}\sim c\circ u$. The fiber-cofiber sequence
\begin{equation*}
   i^*\Pc \xrightarrow{\varepsilon} i^*\Pc \otimes \Lc^{\vee}_{\Vu(p^+)}[2]\xrightarrow{d[1]} i^*i_*i^*\Pc[1] \xrightarrow{c[1]}
\end{equation*}
implies, in particular, the existence of a $2$-cell $\tau$ providing an homotopy $c \circ d\sim 0$. Consider the diagram
\begin{equation*}
    \begindc{\commdiag}[20]
    \obj(-60,10)[1]{$i^*\Pc \otimes \Lc_{\Vu(p^+)}^{\vee}[1]$}
    \obj(0,10)[2]{$i^*i_*i^*\Pc$}
    \obj(60,10)[3]{$i^*\Pc$}
    \obj(-60,-10)[4]{$i^*\Pc \otimes \Lc_{\Vu(p^+)}^{\vee}[1]$}
    \obj(0,-10)[5]{$i^*\Pc \otimes \Lc_{\Vu(p^+)}^{\vee}[1]\oplus i^*\Pc$}
    \obj(60,-10)[6]{$i^*\Pc,$}
    \mor{1}{2}{$d$}
    \mor{2}{3}{$c$}
    \mor{4}{5}{$\begin{bmatrix} 1 \\ 0 \end{bmatrix}$}[\atright,\solidarrow]
    \mor{5}{6}{$\begin{bmatrix} 0 & 1 \end{bmatrix}$}[\atright,\solidarrow]
    \mor{4}{1}{$1$}
    \mor{5}{2}{$\begin{bmatrix} d & u \end{bmatrix}$}
    \mor{6}{3}{$1$}
    \enddc
\end{equation*}
where horizontal lines are fiber-cofiber sequences.
The $2$ cells $\sigma$ and $\tau$ provide an homotopy $c\circ \begin{bmatrix} d & u \end{bmatrix}\sim 1\circ \begin{bmatrix} 0 & 1 \end{bmatrix}$. The existence of an homotopy making the square on the left commutative is clear. Since the left and right vertical morphisms are equivalences, so is the one in the middle. In particular, we obtain that $\varepsilon \sim 0: i^*\Pc \rightarrow \Pc \otimes \Lc_{\Vu(p^+)}^{\vee}[2]$. We have proven that $\ising{X,p}{i}\subseteq \sing{X,p^{+}}^{\varepsilon-tors}$. 

For the converse inclusion, we will show by induction that $fiber(\Mcc \xrightarrow{\varepsilon^n} \Mcc\otimes \Lc^{\vee,\otimes n}_{\Vu(p^+)}[2n]) \in \ising{X,p}{i}$ for every $n$ and for every $\Mcc \in \sing{X,p^+}$. 
\begin{itemize}
    \item The case $n=1$ is obvious, after the definition of $\varepsilon$.
    \item Let $n\geq 1$. The octahedron axiom applied to $\varepsilon^n\sim \varepsilon^{n-1}\circ \varepsilon$ provides us with the following fiber-cofiber sequence:
    \begin{equation*}
        fiber(\Mcc \xrightarrow{\varepsilon} \Mcc \otimes \Lc^{\vee}_{\Vu(p^+)}[2])\rightarrow  fiber(\Mcc \xrightarrow{\varepsilon^n} \Mcc \otimes (\Lc^{\vee}_{\Vu(p^+)})^{\otimes n}[2n])
        \end{equation*}
        \begin{equation*}
        \rightarrow  fiber(\Mcc \otimes \Lc^{\vee}_{\Vu(p^+)}[2] \xrightarrow{\varepsilon^{n-1}} \Mcc \otimes (\Lc^{\vee}_{\Vu(p^+)})^{\otimes n}[2n]).
    \end{equation*}
    Since 
    \begin{equation*}
        fiber(\Mcc \otimes \Lc^{\vee}_{\Vu(p^+)}[2] \xrightarrow{\varepsilon^{n-1}} \Mcc \otimes (\Lc^{\vee}_{\Vu(p^+)})^{\otimes n}[2n])\simeq \Cc \otimes \Lc_{\Vu(p^+)}^{\vee}[2],
    \end{equation*}
    where $\Cc \simeq fiber(\Mcc \xrightarrow{\varepsilon^{n-1}} \Mcc \otimes (\Lc^{\vee}_{\Vu(p^+)})^{\otimes n-1}[2n-2])$, and we have a fiber cofiber sequence
    \begin{equation*}
        i^*i_*\Cc \rightarrow \Cc \xrightarrow{\varepsilon}\Cc \otimes \Lc_{\Vu(p^+)}^{\vee}[2],
    \end{equation*}
    we conclude by thickness of $\ising{X,p}{i}$ that $fiber(\varepsilon^n)\in \ising{X,p}{i}$.
\end{itemize}
Let $\Mcc \in \sing{X,p^+}^{\varepsilon-tors}$. Let $n\geq0$ such that $\varepsilon^n\sim 0$. Then 
\begin{equation*}
    \Mcc \oplus \Mcc \otimes (\Lc_{\Vu(p^+)}^{\vee})^{\otimes n}[2n-1]=fiber(\varepsilon^n)
\end{equation*}
lies in $\ising{X,p}{i}$. Being a retract of a object in $\ising{X,p}{i}$, $\Mcc$ lies in $\ising{X,p}{i}$ as well.
\end{proof}
\begin{cor}
The $\oo$-functor $\singpf{Y}$ enhances the assignment 
\begin{equation*}
    (X,p^+)\rightarrow \ising{X,p}{i}
\end{equation*} 
with the structure of an $\oo$-functor.
\end{cor}
\begin{proof}
This follows immediately from the fact that $\varepsilon$ is functorial with respect to $(X,p^+)$ and from the previous proposition.
\end{proof}
\begin{cor}
There exists an $\oo$-functor 
\begin{equation*}
    \tlgmpop{Y}\rightarrow \textup{Fun}\bigl ( \Delta^1, \dgcatm{S} \bigr )
\end{equation*}
which at the level of objects is defined as 
\begin{equation*}
    (X,p^+)\mapsto \ising{X,p}{i}\hookrightarrow \sing{X,p^+}.
\end{equation*}
By composing it with the $\oo$-functor $\textup{Fun}\bigl ( \Delta^1, \dgcatm{S} \bigr ) \rightarrow \dgcatm{S}$ induced by "taking quotients", we obtain an $\oo$-functorial enhancement 
\begin{equation*}
    \singcf{Y}: \tlgmpop{Y} \rightarrow \dgcatm{S}
\end{equation*}
of the assignment 
\begin{equation*}
    (X,p^+)\mapsto \singcf{Y}(X,p^+)=\singc{\Vu(p),s}=\frac{\sing{X,p^+}}{\ising{X,p}{i}}.
\end{equation*}
\end{cor}
\begin{comment}
\red{Can we enhance it with a lax monoidal structure?}
\end{comment}

Let us now focus on the other side of the equivalence provided by Theorem \ref{reduction of codimension for sing}, where everything we said before applies, \textit{mutatis mutandis}. There exists a lax monoidal $\oo$-functor
\begin{equation*}
    \singprojfm{\Ec_Y^{+,\vee}}:\tlgmprojmop{\Ec_Y^{+,\vee}}\rightarrow \dgcatmo{S}
\end{equation*}
which, at the level of objects, is defined by
\begin{equation*}
    (Z,w)\mapsto \singprojf{\Ec_Y^{+,\vee}}(Z,w)=Ker\bigl ( g_*: \sing{\Vu(w)}\rightarrow \sing{Z}\bigr ),
\end{equation*}
where $g$ denotes the closed embedding of the zero locus of $w\in \Hu^0(Z,\Oc(1))$ inside $Z$. Similarly, there exists an analogous lax monoidal $\oo$-functor
\begin{equation*}
    \singprojfm{\Ec_Y^{\vee}}:\tlgmprojmop{\Ec_Y^{\vee}}\rightarrow \dgcatmo{S}.
\end{equation*}
We will consider their compositions with the symmetric monodail functors $\Xi^{\boxplus}$ and $Q^{\boxplus}$ (see Remark \ref{functor Q} and the discussion thereafter):
\begin{equation*}
    \singprojfm{\Ec_Y^{+,\vee}}\circ \Xi^{\boxplus} :\tlgmpmop{Y} \rightarrow \dgcatmo{S},
\end{equation*}
\begin{equation*}
    \singprojfm{\Ec_Y^{\vee}}\circ \Xi^{\boxplus}\circ Q^{\boxplus}: \tlgmpmop{Y}\rightarrow \dgcatmo{S}.
\end{equation*}
\begin{prop}
Let $(X,p^+)$ be an LG model over $(Y,\Ec_Y^+)$ and let $r:\Vu(W_p)\rightarrow \Vu(W_{p^+})$ denote the canonical lci closed embedding. The adjunction
\begin{equation*}
    r_*:\sing{\Pbb(\Ec_X^{\vee}),W_p}\rightleftarrows \sing{\Pbb(\Ec_X^{+,\vee}),W_{p^+}}:r^!
\end{equation*}
is functorial in $(X,p^+)$, i.e. we have two natural transformations
\begin{equation*}
\singprojf{\Ec_Y^{\vee}}\circ \Xi \circ Q\rightarrow \singprojf{\Ec_Y^{+,\vee}}\circ \Xi, \hspace{0.3cm} \singprojf{\Ec_Y^{+,\vee}}\circ \Xi \rightarrow \singprojf{\Ec_Y^{\vee}}\circ \Xi \circ Q 
\end{equation*}
that, for a fixed LG model $(X,p^+)$ over $(Y,\Ec_Y^+)$, determine the above mentioned adjunction.

In particular, the BW-operator is functorial in $(X,p^+)$.
\end{prop}
\begin{proof}
This is clear when we restrict to pairs $(X,p^+)$ with $p^+$ regular. Indeed, we have proved the functoriality on $(X,p^+)$ for
\begin{equation*}
    i^*\bigl ((-)\otimes \Lc_{\Vu(p)}  \bigr )[-1]:\sing{X,p^+}\rightleftarrows \sing{X,p}:i_*\bigl ((-)\otimes \Lc_{\Vu(p)}^{\vee}  \bigr )[1]
\end{equation*}
and the equivalences $\Upsilon_{(X,p)}$ and $\Upsilon_{(X,p^+)}$ are likewise functorial in $(X,p^+)$. However, this argument can't be used in our setup as we have not proved that $\Upsilon_{(X,p)}$ and $\Upsilon_{(X,p^+)}$ are always equivalences, regardless of the fact that the section is regular or not. 

Nevertheless, we can mimic the proof of Proposition \ref{functoriality (i^*,i_*)}: using strict models, one can show that $r^!$ is functorial in $(X,p^+)$. Indeed, the same argument applies \emph{mutatis mutandis}. The only difference is that the $r^!$ is modelled a pseudo-functor of the form
\begin{equation*}
    \cohs{B,M^{\vee}_B,w}\rightarrow \cohs{B,(M_B \oplus L_B)^{\vee},(w,t)}
\end{equation*}
\begin{equation*}
    F\mapsto F\otimes K(L_B^{\vee},t)\otimes_B O_B[-1],
\end{equation*}
where $O_B$ is a line bundle (in fact, the restriction of $\Oc(1)\otimes \Lc^{\vee}_{\Vu(W_p)}$) such that, for any $B\rightarrow C$, there is a canonical equivalence $O_C\simeq O_B\otimes_BC$. This doesn't affect the validity of the proof.
\end{proof}
We will now focus on the study of $\chi$-torsion objects inside $\sing{\Pbb(\Ec_X^{+,\vee}),W_{p^+}}$.
\begin{defn}
Let $(X,p^+)$ be an LG model over $(Y,\Ec_Y^+)$. The dg category of $\chi$-torsion objects in $\sing{\Pbb(\Ec_X^{+,\vee}),W_{p^+}}$ is the full sub dg category 
\begin{equation*}
    \sing{\Pbb(\Ec_X^{+,\vee}),W_{p^+}}^{\chi-tors}:=
\end{equation*}
\begin{equation*}
    \Bigl \{\Cc \in \sing{\Pbb(\Ec_X^{+,\vee}),W_{p^+}}: \exists n\geq0, \chi^n\sim 0: \Cc \rightarrow \Cc \otimes \Oc_{\Vu(W_{p^+})}(n)\otimes \Lc_{\Vu(W_{p^+})}^{\vee,\otimes n} \Bigr \}.
\end{equation*}
\end{defn}
\begin{prop}\label{BW-torsion = image r_*}
Let $(X,p^+)$ be an LG model over $(Y,\Ec_Y^+)$. Let $\rsing{\Pbb(\Ec_X^{\vee}),W_p}{r}$ denote the thick full sub dg category of $\sing{\Pbb(\Ec_X^{+,\vee}),W_{p^+}}$ spanned by the image of the dg functor
\begin{equation*}
    r_*:\sing{\Pbb(\Ec_X^{\vee}),W_p}\rightarrow \sing{\Pbb(\Ec_X^{+,\vee}),W_{p^+}}.
\end{equation*}
Then,
\begin{equation*}
    \sing{\Pbb(\Ec_X^{+,\vee}),W_{p^+}}^{\chi-tors}\simeq \rsing{\Pbb(\Ec_X^{\vee}),W_p}{r}.
\end{equation*}
\end{prop}
\begin{proof}
The proof of Proposition \ref{E-torsion = image pullback} works in this case too. Assume that $\Pc \in \sing{\Pbb(\Ec_X^{\vee}),W_p}$. The counit morphism $r_*r^!r_*\Pc \xrightarrow{c} r_*\Pc$ admits a section $r_*\Pc \xrightarrow{u}r_*r^!r_*\Pc$. Therefore, as in the proof of Proposition \ref{E-torsion = image pullback}, we obtain that
\begin{equation*}
    r_*r^!r_*\Pc\simeq r_*\Pc\oplus r_*\Pc\otimes \Oc(1)\otimes \Lc^{\vee}[-1].
\end{equation*}
In particular, this means that $\chi \sim 0: r_*\Pc \rightarrow r_*\Pc \otimes \Oc(1)\otimes \Lc^{\vee}$ and therefore $r_*\Pc$ lies in  $\sing{\Pbb(\Ec_X^{+,\vee}),W_{p^+}}^{\chi-tors}$.

For the converse, we shall first show that $fiber(\Mcc \xrightarrow{\chi^n}\Mcc \otimes \Oc(n)\otimes \Lc^{\vee,\otimes n}) \in \rsing{\Pbb(\Ec_X^{\vee}),W_p)}{r}$ for any $n\geq 0$. We proceed by induction on $n$:
\begin{itemize}
    \item the case $n=1$ follows from the definition of $\chi$.
    \item if $n\geq 1$, the claim follows from the fiber-cofiber sequence associated to $\chi^n\sim \chi^{n-1}\circ \chi$ .
\end{itemize}
Then, if $\chi^n\sim 0:\Mcc \rightarrow \Mcc \otimes \Oc(n)\otimes \Lc^{\vee,\otimes n}$, we see that $\Mcc$ is a retract of an object in $\rsing{\Pbb(\Ec_X^{\vee}),W_p}{r}$.
\end{proof}
\begin{cor}
The $\oo$-functor $\singprojf{\Ec_Y^{\vee}}$ enhances the assignment 
\begin{equation*}
    (X,p^+)\rightarrow \rsing{\Pbb(\Ec_X^{\vee}),W_p}{r}
\end{equation*} 
with the structure of an $\oo$-functor.
\end{cor}
\begin{proof}
This follows immediately from the fact that $\chi$ is functorial with respect to $(X,p^+)$ and from the previous proposition.
\end{proof}
\begin{cor}
There exists an $\oo$-functor 
\begin{equation*}
    \tlgmpop{Y}\rightarrow \textup{Fun}\bigl ( \Delta^1, \dgcatm{S} \bigr )
\end{equation*}
which, at the level of objects, is defined as 
\begin{equation*}
    (X,p^+)\mapsto \rsing{\Pbb(\Ec_X^{\vee}),W_p}{r}\hookrightarrow \sing{\Pbb(\Ec_X^{+,\vee}),W_{p^+}}.
\end{equation*}
By composing it with the quotient $\oo$-functor $\textup{Fun}\bigl ( \Delta^1, \dgcatm{S} \bigr ) \rightarrow \dgcatm{S}$, we obtain an $\oo$-enhancement 
\begin{equation*}
    \singprojcf{\Ec_Y^{+,\vee}}: \tlgmpop{Y} \rightarrow \dgcatm{S}
\end{equation*}
of the assignment 
\begin{equation*}
    (X,p^+)\mapsto \frac{\sing{\Pbb(\Ec_X^{+,\vee}),W_{p^+}}}{\rsing{\Pbb(\Ec_X^{\vee})}{r}}.
\end{equation*}
\end{cor}
\begin{cor}
Let $(X,p^+)$ be an LG model over $(Y,\Ec_Y^+)$ and assume that $p^+$ is a regular section. Then we have equivalences
\begin{equation*}
    \ising{X,p}{i}\simeq \sing{X,p^+}^{\varepsilon-tors}\simeq \sing{\Pbb(\Ec_X^{+,\vee}),W_{p^+}}^{\chi-tors}\simeq \rsing{\Pbb(\Ec_X^{\vee}),W_p}{r}.
\end{equation*}
\end{cor}
\begin{proof}
The first equivalence follows from Proposition \ref{E-torsion = image pullback}, the second from Corollary \ref{E-operator = BW-operator} and the last one from Proposition \ref{BW-torsion = image r_*}.
\end{proof}
\begin{cor}
Let $(X,p^+)$ be an LG model over $(Y,\Ec_Y^+)$ and assume that $p^+$ is a regular section. Then
\begin{equation*}
    \singc{\Vu(p),s}\simeq \frac{\sing{\Pbb(\Ec_X^{+,\vee}),W_{p^+}}}{\rsing{\Pbb(\Ec_X^{\vee}),W_p}{r}}.
\end{equation*}
\end{cor}
Our next aim is to better understand the dg category $\rsing{\Pbb(\Ec_X^{\vee}),W_p}{r}$. Indeed, at least for pairs $(X,p^+)$ where $p^+$ is regular, the results of the previous section imply that we have an exact sequence of dg categories
\begin{equation*}
    \rsing{\Pbb(\Ec_X^{\vee}),W_p}{r} \hookrightarrow \sing{\Pbb(\Ec_X^{+,\vee}),W_{p^+}}\rightarrow \singc{\Vu(p),s},
\end{equation*}
which is a tool to study the certain (co)homologies of $\singc{\Vu(p),s}$.
\begin{context}
From now on we shall assume that $(X,p^+)$ is an LG model over $(Y,\Ec_Y^+)$ with $X$, $Y$ and $p^+$ regular.
\end{context}
\begin{rmk}
Under our standing assumptions, it is immediate to observe that
\begin{equation*}
    \sing{X,p^+}\simeq \sing{\Vu(p^+)}, \hspace{0.5cm} \sing{X,p}\simeq \sing{\Vu(p)},
\end{equation*}
\begin{equation*}
    \sing{\Pbb(\Ec_X^{+,\vee}),W_{p^+}}\simeq \sing{\Vu(W_{p^+})}, \hspace{0.5cm} \sing{\Pbb(\Ec_X^{\vee}),W_{p}}\simeq \sing{\Vu(W_{p})}.
\end{equation*}
In particular, the exact sequence of dg categories above looks as follows
\begin{equation*}
     \rsing{\Vu(W_p)}{r} \hookrightarrow \sing{\Vu(W_{p^+})}\rightarrow \singc{\Vu(p),s}.
\end{equation*}
\end{rmk}
First of all, notice the existence of the following diagram of dg categories:
\begin{equation*}
    \begindc{\commdiag}[20]
    \obj(-60,10)[1]{$\rperf{\Vu(W_p)}{r}$}
    \obj(0,10)[2]{$\rcoh{\Vu(W_p)}{r}$}
    \obj(60,10)[3]{$\rsing{\Vu(W_p)}{r}$}
    \obj(-60,-10)[4]{$\perfsupp{\Vu(W_{p^+})}{\Vu(W_p)}$}
    \obj(0,-10)[5]{$\cohbsupp{\Vu(W_{p^+})}{\Vu(W_p)}$}
    \obj(60,-10)[6]{$\singsupp{\Vu(W_{p^+})}{\Vu(W_p)}.$}
    \mor{1}{2}{$$}[\atright,\injectionarrow]
    \mor{2}{3}{$$}
    \mor{1}{4}{$$}[\atright,\injectionarrow]
    \mor{2}{5}{$$}[\atright,\injectionarrow]
    \mor{3}{6}{$$}[\atright,\injectionarrow]
    \mor{4}{5}{$$}[\atright,\injectionarrow]
    \mor{5}{6}{$$}
    \enddc
\end{equation*}

The dg categories on the bottom line are the full sub dg categories of $\perf{\Vu(W_{p^+})}$ (resp. $\cohb{\Vu(W_{p^+})}$, resp. $\sing{\Vu(W_{p^+})}$) spanned by those objects with support contained in the closed subset $\Vu(W_p)$.

\begin{rmk}
Notice that the BW-operator is well defined on $\cohb{\Vu(p^+)}$ and on $\perf{\Vu(p^+)}$. In particular, we can consider $\chi$-torsion objects in $\perf{\Vu(W_{p^+})}$ and in $\cohb{\Vu(W_{p^+})}$:
\begin{equation*}
    \perf{\Vu(W_{p^+})}^{\chi-tors}=\bigl \{ \Cc \in \perf{\Vu(W_{p^+})}: \exists n\geq0, \chi^n\sim 0:\Cc\rightarrow \Cc \otimes \Oc(n) \otimes \Lc_{\Vu(W_{p^+})}^{\vee,\otimes n} \bigr \},
\end{equation*}
\begin{equation*}
    \cohb{\Vu(W_{p^+})}^{\chi-tors}=\bigl \{ \Cc \in \cohb{\Vu(W_{p^+})}: \exists n\geq0, \chi^n\sim 0:\Cc\rightarrow \Cc \otimes \Oc(n)\otimes \Lc_{\Vu(W_{p^+})}^{\vee,\otimes n} \bigr \}.
\end{equation*}
\end{rmk}
\begin{prop}
The following equivalences hold:
\begin{enumerate}
    \item $\rperf{\Vu(W_p)}{r}\simeq \perf{\Vu(W_{p^+})}^{\chi-tors}\simeq \perfsupp{\Vu(W_{p^+})}{\Vu(W_p)}$;
    \item $\rcoh{\Vu(W_p)}{r}\simeq \cohb{\Vu(W_{p^+})}^{\chi-tors}\simeq \cohbsupp{\Vu(W_{p^+})}{\Vu(W_p)}$;
    \item $\rsing{\Vu(W_p)}{r}\simeq \sing{\Vu(W_{p^+})}^{\chi-tors}\simeq \singsupp{\Vu(W_{p^+})}{\Vu(W_p)}$.
\end{enumerate}
\end{prop}
\begin{proof}
The first equivalences in (1), (2) and (3) can be proved as Proposition \ref{E-torsion = image pullback}. The other equivalences follow from the observation that $\chi:\Cc \rightarrow \Cc \otimes \Oc(1) \otimes \Lc_{\Vu(W_{p^+})}^{\vee}$ is homotopic to the map induced by $T$, the global section of $\Oc(1)\otimes \Lc_{\Vu(W_{p^+})}^{\vee}$ which defines the closed embedding $\Vu(W_p)\rightarrow \Vu(W_{p^+})$. 
Indeed, one can easily see that there is an equivalence
\begin{equation*}
    r_*r^*(\Cc)\simeq cofib \bigl ( \Cc \otimes \Oc(-1)\otimes \Lc_{\Vu(W_{p^+})} \xrightarrow{T}\Cc \bigr ).
\end{equation*}
Using the equivalences $r^!(\Cc)\simeq r^*\Cc \otimes \omega_{\Vu(W_p)/\Vu(W_{p})}\simeq r^*\Cc \otimes \Oc(1)\otimes \Lc_{\Vu(W_{p})}^{\vee}[-1]$, we get that
\begin{equation*}
    r_*r^!(\Cc)\simeq r_*r^*\bigl (\Cc \otimes \Oc(1)\otimes \Lc^{\vee}_{\Vu(W_{p^+})} \bigr )[-1]\simeq fib \bigl ( \Cc \xrightarrow{T} \Cc \otimes \Oc(1)\otimes \Lc^{\vee}_{\Vu(W_{p^+})} \bigr ).
\end{equation*}
Then, $\Cc$ is $\chi$-torsion object if and only if $T^n:\Cc \rightarrow \Cc \otimes \Oc(n) \otimes \Lc_{\Vu(W_{p^+})}^{\vee,\otimes n}$ is homotopic to zero for some $n>>0$, i.e. the support of $\Cc$ is contained in $\Vu(W_p)$.
\end{proof}
\begin{cor}
The following is an exact sequence of dg categories:
\begin{equation*}
    \rperf{\Vu(W_p)}{r}\hookrightarrow \rcoh{\Vu(W_p)}{r} \rightarrow \rsing{\Vu(W_p)}{r}.
\end{equation*}
\end{cor}
\begin{proof}
This is clear from the previous proposition, as
\begin{equation*}
    \perfsupp{\Vu(W_{p^+})}{\Vu(W_p)}\hookrightarrow \cohbsupp{\Vu(W_{p^+})}{\Vu(W_p)}\rightarrow \singsupp{\Vu(W_{p^+})}{\Vu(W_p)}
\end{equation*}
is an exact sequence of dg categories.
\end{proof}
In particular, we obtain the following result, which can be thought as the main result of this note:
\begin{cor}\label{main equivalence}
The equivalence $\Upsilon_{(X,p^+)}$ induces an equivalence of exact sequences in $\dgcatm{S}$:
\begin{equation*}
    \begindc{\commdiag}[20]
    \obj(-60,10)[1]{$\ising{\Vu(p)}{i}$}
    \obj(0,10)[2]{$\sing{\Vu(p^+)}$}
    \obj(60,10)[3]{$\singc{\Vu(p),s}$}
    \obj(-60,-10)[4]{$\singsupp{\Vu(W_{p^+})}{\Vu(W_p)}$}
    \obj(0,-10)[5]{$\sing{\Vu(W_{p^+})}$}
    \obj(60,-10)[6]{$\sing{\Vu(W_{p^+})-\Vu(W_p)}.$}
    \mor{1}{2}{$$}[\atright,\injectionarrow]
    \mor{2}{3}{$$}
    \mor{1}{4}{$$}
    \mor{2}{5}{$$}
    \mor{3}{6}{$$}
    \mor{4}{5}{$$}[\atright,\injectionarrow]
    \mor{5}{6}{$$}
    \enddc
\end{equation*}
\end{cor}

\begin{rmk}
Notice that, when $\Vu(p)$ is a regular scheme, we have the implication
\begin{equation*}
    \sing{\Vu(W_p)}\simeq 0 \Rightarrow \singsupp{\Vu(W_{p^+})}{\Vu(W_p)}\simeq 0.
\end{equation*}
It is easy to see that this phenomenon is false in general. Given a closed embedding $Z\rightarrow Y$ with $Z$ regular, it is not always true that $\singsupp{Y}{Z}\simeq 0$. For example, consider the closed embedding of the singular point of a nodal curve inside the curve.
\end{rmk}

\section{Applications}\label{applications}
\subsection{\texorpdfstring{$\ell$}{l}-adic realization of \texorpdfstring{$\singc{\Vu(p),s}$}{Singcoh}}

Let $(X,p^+)$ be an LG model over $(Y,\Ec_Y^+)$.
Assume that both $X$ and $p^+$ are regular. One immediate consequence of the identification
\begin{equation*}
    \singc{\Vu(p),s}\simeq \sing{\Vc,W_{p^+|\Vc}}\simeq \sing{\Uc},
\end{equation*}
where $\Vc:=\Pbb(\Ec_X^{+,\vee})-\Pbb(\Ec^{\vee}_X)$, $W_{p^+|\Vc}$ denotes the restiction of the global section $W_{p^+}\in \Hu^0(\Pbb(\Ec_X^{+,\vee}),\Oc(1))$ to $\Vc$ and $\mathcal{U}=\Vu(W_{p^+|\Vc})=\Vu(W_{p^+})-\Vu(W_p)$, is the computation of the $\ell$-adic cohomology of the dg category of coherent matrix factorizations. 
\begin{rmk}
Using the isomorphisms
\begin{equation*}
    \Pbb(\Ec_X^{\vee})\simeq \Pbb(\Ec_X^{\vee}\otimes \Lc_X), \hspace{0.5cm} \Pbb(\Ec_X^{+,\vee})\simeq \Pbb(\Ec_X^{\vee}\otimes \Lc_X \oplus \Oc_X),
\end{equation*}
we can see that 
\begin{equation*}
    \Vc \simeq \V(\Ec_X^{\vee}\otimes \Lc_X)=Spec_X(Sym_{\Oc_X}(\Ec_X\otimes \Lc_X^{\vee})).
\end{equation*}
\end{rmk}
Let $\ell$ be an invertible prime number in $S$. Recall (see \cite{brtv18} and \cite{p20}) that, for any scheme $Z$ of finite type over $S$, there is a lax monoidal $\oo$-functor
\begin{equation*}
    \Rlv_Z:\dgcatmo{Z} \rightarrow \shvl(Z)^{\otimes},
\end{equation*}
which sends a dg category $\Cc$ to an $\ell$-adic sheaf $\Rlv(\Cc)$, which we refer to as the $\ell$-adic cohomology of $\Cc$. This $\oo$-functor is compatible with the usual $\ell$-adic realization of schemes by means of the following equivalence
\begin{equation*}
    \Rlv_Z(\perf{Y})\simeq f_*\Ql{,Y}(\beta),
\end{equation*}
where $f:Y\rightarrow Z$ and $\Ql{,Y}(\beta)\simeq \bigoplus_{n\in \Z}\Ql{,Y}(n)[2n]$.

Also recall the notion of monodromy-invariant vanishing cycles (see \cite{p20}): let $(Z,\Mcc_Z)$ be a scheme of finite type over $S$ together with a line bundle. Let $s \in \Hu^0(Z,\Mcc_Z)$ be a regular section and consider the diagram
\begin{equation*}
    \begindc{\commdiag}[20]
    \obj(-50,20)[1]{$\Vu(s)$}
    \obj(0,20)[2]{$Z$}
    \obj(50,20)[3]{$Z-\Vu(s)$}
    \obj(-50,-20)[4]{$Z$}
    \obj(0,-20)[5]{$\V(\Mcc_Z)$}
    \obj(50,-20)[6]{$\V(\Mcc_Z)-Z.$}
    \mor{1}{2}{$i$}
    \mor{3}{2}{$j$}[\atright,\solidarrow]
    \mor{1}{4}{$s_0$}
    \mor{2}{5}{$s$}
    \mor{3}{6}{$s_U$}
    \mor{4}{5}{$i_0=0$}[\atright,\solidarrow]
    \mor{6}{5}{$j_0$}
    \enddc
\end{equation*}
We define the $\ell$-adic sheaf of monodromy invariant vanishing cycles of $(Z,s)$ as
\begin{equation*}
    \Phi_{(Z,s)}^{\textup{mi}}(\Ql{}(\beta)):=cofib \bigl ( i^*s^*j_{0*}\Ql{,\V(\Mcc_Z)-Z}(\beta)\rightarrow i^*j_*s_{U}^*\Ql{,\V(\Mcc_Z)-Z}(\beta)\bigr ).
\end{equation*}
\begin{thm}\label{l-adic cohomology mf coh}
Let $(X,p^+)$ be an LG model over $(Y,\Ec_Y^+)$ such that $X$ is a regular scheme and $p^+$ is a regular section. Then
\begin{equation*}
   \Rlv_{Y}(\singc{\Vu(p),s})\simeq \Rlv_{Y}(\mfc{\Vu(p)}{s}\simeq q_*\Phi_{(\Vc,W_{p^+|\Vc})}^{\textup{mi}}(\Ql{}(\beta))[-1],
\end{equation*}
where $q:\Uc \rightarrow Y$ is the canonical morphism.
\end{thm}
\begin{proof}
This follows immediately from the discussion above and from \cite[Theorem 5.2.2]{p20}
\end{proof}
\begin{rmk}
The equivalence of the theorem above stands at the motivic level, i.e. before taking the $\ell$-adic realization:
\begin{equation*}
     \Mv_{Y,\Q}(\singc{\Vu(p),s})\simeq \Mv_{Y,\Q}(\mfc{\Vu(p)}{s}) \simeq  q_*\Phi_{(\Vc,W_{p^+|\Vc})}^{\textup{mi, mot}}(\bu_{\Q})[-1].
\end{equation*}
More precisely, the equivalence above lives in $\md_{\bu_{\Q}}(\textup{\textbf{SH}}_Y)$, the $\oo$-category of modules over the spectrum of homotopy invariant, non connective rational K theory. See \cite{brtv18} and \cite{p20} for notation.
\end{rmk}
The following example seems particularly interesting.
\begin{ex}
Let $S$ be an excellent, strictly henselian trait and let $X$ be a projective $S$-scheme. Then there exist an integer $N$ and homogenous polynomials $p \in \Hu^0(Y,\oplus_{k=1}^m \Oc_Y(d_k))$, where $Y=\proj{N}{S}$ and $d_k\geq 0$, $k=1,\dots,m$, such that $X=V(p)$. We assume that $p$ is a regular section of $\bigoplus_{k=1}^m\Oc_Y(d_k)$. Fix an uniformizer $\pi$ of $S$. We shall label $\pi_Y$ the induced regular function on $\proj{N}{S}$. Setting $\Ec_Y=\oplus_{k=1}^m \Oc_Y(d_k)$, $\Ec_Y^+=\Ec_Y\oplus \Oc_Y$ and $p^+=(p,\pi_Y)$, it is immediate to see that $(Y,p^+)$ defines a LG model over $(Y,\Ec_Y^+)$ such that $\Vu(p)=X$ and $\Vu(p^+)=X_{\sigma}$, the special fiber of $X$ over $S$, i.e. the zero locus of the pullback $\pi_X$ of $\pi$ along $X\rightarrow S$. Then all the previous results apply and we find that
\begin{equation*}
    \mfc{X}{\pi_X}\simeq \sing{\Uc},
\end{equation*}
where $\Uc=\Vu(W_{p^+|\Vc})$ and $\mfc{X}{\pi_X}$ denotes the dg category of coherent matrix factorizations of $(X,\pi_X)$ (see \cite{ep15}). In particular, we get
\begin{equation*}
    \Rlv_Y(\mfc{X,\pi_X})\simeq q_*\Phi_{(\Vc,W_{p^+|\Vc})}^{\textup{mi}}(\Ql{}(\beta))[-1].
\end{equation*}
An interesting $\ell$-adic sheaf attached to a an $S$-scheme is that of vanishing cycles. It is known (see \cite[Theorem 4.39]{brtv18}) that when $X$ is regular one can express the (inertia invariant) vanishing cohomology of $X\rightarrow S$ by means of the singularity category of the special fiber, which agrees with $\mfc{X}{\pi_X}$:
\begin{equation*}
    i^*_{\sigma}\Rlv_S(\mfc{X}{\pi_X}) \simeq i^*_{\sigma}\Rlv_S(\sing{X_{\sigma}})\simeq\bigl (p_{\sigma*}\nu(\Ql{,X}(\beta))\bigr )^{\textup{h}I}[-1],
\end{equation*}
where $I$ denotes the inertia group of $S$, $i_{\sigma}:\sigma \rightarrow S$ is the closed point, $p_{\sigma}:X_{\sigma}=X\times_S \sigma \rightarrow \sigma$ the projection and $\nu(\Ql{,X}(\beta))$ is the $\ell$-adic sheaf of vanishing cycles of $\Ql{,X}(\beta)$. Recall that $\Ql{,X}(\beta)$ is the commutative algebra object $Sym(\Ql{,X}\cdot \beta)[\beta^{-1}]$ with $\beta$ living in bidegree $(1,2)$, i.e. $\Ql{,X}(\beta)=\bigoplus_{n\in \Z}\Ql{,X}(n)[2n]$. 

This implies, in particular, that there is an equivalence
\begin{equation*}
    \Phi_{(\Vc,W_{p^+})}^{\textup{mi}}(\Ql{}(\beta))\simeq \bigl ( \nu(\Ql{,X}(\beta))\bigr )^{\textup{h}I}
\end{equation*}
whenever $X$ is regular. It would be interesting to understand whether there exists an $\ell$-adic sheaf $\Fc$ on $X$, which agrees with $\Ql{,X}$ if $X$ is regular, such that the equivalence 
\begin{equation*}
    \Phi_{(\Vc,W_{p^+})}^{\textup{mi}}(\Ql{}(\beta))\simeq \bigl ( \nu(\Fc(\beta)) \bigr )^{\textup{h}I}
\end{equation*}
holds without the regularity assumptions.
\end{ex}
\subsection{Hochschild and periodic cyclic homologies of \texorpdfstring{$\singc{\Vu(p),s}$}{Singcoh}}
The connection between categories of singularities and vanishing cycles is well known and predates the above mentioned theorem of Blanc, Robalo, Toën and Vezzosi. For example, Hochschild (co)homology of matrix factorizations has been computed by Efimov (\cite{ef18}), Dyckerhoff (see \cite{dy11}) Lin-Pomerlano (see \cite{lp13}), Preygel (see \cite{pr11}), Segal (see \cite{se13}). However, the author is not aware of any computation of these invariants for $\mfc{X}{f}$ (or $\mf{X}{f}$) when $X$ is not supposed regular.

Let us explain these results as stated in \cite{ef18}. Let $\Cc$ be a $\Z/2\Z$-graded dg category. Its \textit{Hochschild homology} is defined as
\begin{equation*}
    \hh{\Cc}:=\Hu^{-\bullet}(id_{\Cc}\otimes_{\Cc \otimes \Cc^{\textup{op}}}id_{\Cc}),
\end{equation*}
where $id_{\Cc}$ denotes the identity $\Cc \otimes \Cc^{\textup{op}}$-bimodule. It is computed by the bar complex $\hoch{\Cc}{b}$ as defined e.g. in \cite[\S3]{ef18}.

The notions of \emph{mixed complex} and \emph{u-connection on a mixed complex} are crucial to understand Efimov's results. Let us recall these notions for the reader's convenience.
\begin{defn} Let $k$ be a field.
\begin{itemize}
    \item A \emph{mixed complex} is a triple $(C,b,B)$, where $C$ is a $\Z/2\Z$-graded $k$-vector space and $b$ and $B$ are two odd differentials on $C$ such that
    \begin{equation*}
        bB+Bb=0.
    \end{equation*}
    \item A morphism $f:(C_1,b_1,B_1)\rightarrow (C_2,b_2,B_2)$ is a graded morphism of $k$-vector spaces which commutes with both differentials.
    \item $f:(C_1,b_1,B_1)\rightarrow (C_2,b_2,B_2)$ is said to be a quasi-isomorphism if $f:(C_1,b_1)\rightarrow (C_2,b_2)$ is a quasi-isomorphism.
\end{itemize}
\end{defn}
\begin{rmk}
Notice that, even if the assignment $(C,b,B)\mapsto (C,B,b)$ defines an endofunctor on the category of mixed complexes, it doesn't preserve quasi-isomorphisms. In other words, the roles played by the two differentials $b$ and $B$ are not symmetric, and $b$ should be thought as the main differential.
\end{rmk}
\begin{ex}
\begin{itemize}
    \item Let $\Cc$ be a $\Z/2\Z$-graded dg category. The \emph{Hochschild complex} $(Hoch(\Cc),b)$, together with the Connes differential $B$, defines a mixed complex. For the precise definition, see \cite[\S3]{ef18}.
    \item Let $\Cc$ be a $\Z/2\Z$-graded curved dg category. The \emph{Hochschild complex of the second kind} $(Hoch^{II}(\Cc),b)$, together with the Connes differential $B$, defines a mixed complex. For the precise definition, see \cite{pp12}, \cite[\S3]{ef18}.
    \item Let $A$ be a smooth $k$ algebra. The \emph{twisted de Rham complex} $(\Omega_{A/k}^{\bullet},-dW\wedge,d_{dR})$ is a mixed complex. 
\end{itemize}
\end{ex}
Let $(C,b,B)$ be a mixed complex. The second differential $B$ is used to define a second complex out of the initial datum. Let $u$ be a formal variable of even degree and set
\begin{equation*}
    C\dsl u \dsr := \varprojlim_{n}C\otimes_k \frac{k[u]}{u^n}, \hspace{0.5cm} C\drl u \drr := \varinjlim_{n} u^{-n}\cdot C\dsl u\dsr.
\end{equation*}
It is immediate that $b+uB$ defines a differential on $C\drl u \drr$ and that every morphism $f:(C_1,b_1,B_1)\rightarrow (C_2,b_2,B_2)$ of mixed complexes induces a morphism of complexes
\begin{equation*}
    f\drl u \drr : \bigl ( C_1\drl u \drr, b_1+uB_1 \bigr )\rightarrow  \bigl ( C_2\drl u \drr, b_2+uB_2\bigr ).
\end{equation*}
Moreover, if $f$ is a quasi-isomorphism of mixed complexes, $f\drl u \drr$ is a quasi-isomorphism too.

Let $B$ denote the Connes differential acting on $\hoch{\Cc}{b}$. The \emph{periodic cyclic homology} of $\Cc$ is defined as
\begin{equation*}
    \hp{\Cc}:= \Hu_{\bullet}\hoch{\Cc \drl u \drr}{b+uB},    
\end{equation*}
where $u$ is a formal variable of even degree.

\begin{defn}
Let $(C,b,B)$ be a mixed complex.
\begin{itemize}
    \item A \emph{u-connection} is a $k$-linear operator
    \begin{equation*}
        \nabla=\frac{d}{du}+M:C\drl u \drr \rightarrow C\drl u \drr,
    \end{equation*}
    where $M$ is a $k\drl u \drr$-linear operator $C\drl u \drr \rightarrow \drl u \drr$, such that
    \begin{equation*}
        [\nabla,b+uB]=\frac{1}{2u}(b+uB).
    \end{equation*}
    \item Let $(C_i,b_i,B_i,\nabla_i)$, $i\in \{1,2\}$ be two mixed complexes with $u$-connections. A morphism of mixed complexes $f:(C_1,b_1,B_1)\rightarrow (C_2,b_2,B_2)$ is \emph{weakly compatible} (resp. \emph{strictly compatible}) \emph{with the u-connections} if $\nabla_2\circ f\drl u \drr-f\drl u \drr \circ \nabla_1$ is homotopic (resp. equal) to zero. The homotopy is part of the datum.
\end{itemize} 
\end{defn}
\begin{rmk}
All the examples of mixed complexes given above can be endowed with $u$-connections in a natural way. In particular, there is a \emph{Getzler-Gauss-Manin u-connection} $\nabla_u^{GM}$ on $\bigl ( Hoch(\Cc),b,B \bigr )$ (see \cite{kkp08} and \cite{sh14}) and a $u$-connection $\nabla_u^{dR}$ on the twisted de Rham complex. We refer the reader to \cite{ef18} for more details.
\end{rmk}
In \textit{loc. cit.} Efimov proves the following
\begin{thm}{\cite[Theorem 1.3, Theorem 1.4]{ef18}}\label{theorem Efimov}
Let $X$ be a separated smooth scheme of finite type over a field $k$ of characteristic zero. Let $W\in \Hu^0(X,\Oc_X)$ be a regular function whose critical locus is contained in the fiber over $0$. There is a chain of quasi-isomorphisms of mixed complexes with $u$-connections between
\begin{equation*}
    \hoch{\mfc{X}{W}}{b,B,\nabla_u^{GM}} 
\end{equation*}
and
\begin{equation*}
    \textup{R}\Gamma \bigl (X, (\Omega_X^{\bullet},-dW,d_{dR},\nabla_u^{dR}) \bigr ).
\end{equation*}
In particular, there is an equivalence between $\Z/2\Z$-graded complexes with connections
\begin{equation*}
    \bigl ( \hp{\mfc{X}{W}},\nabla_u^{GM} \bigr )\simeq \bigl ( \Hu^{\bullet}(\Omega_X^{\bullet}\drl u\drr,-dW +ud_{dR}),\nabla_u^{dR}=\frac{d}{du}+\frac{\Gamma}{u}+\frac{W}{u^2}\bigr ),
\end{equation*}
where $\Gamma_{|\Omega_X^p}=-\frac{p}{2}$, and a quasi-isomorphism of $\Z/2\Z$-graded complexes
\begin{equation*}
    \hh{\mfc{X}{W}}\simeq \Hu^{\bullet}(X,\Omega_X^{\bullet},-dW).
\end{equation*}
\end{thm}
The connection between matrix factorizations and vanishing cycles in this context arises from the combination of Efimov's theorem with the algebraic formula for vanishing cycles, conjectured by  Kontsevich and proved by Sabbah (\cite{sa12})/Sabbah-M.~Saito (\cite{sasa14}):
\begin{thm}{\cite[Theorem 1.1]{sa12}}
For a finite dimensional $\C$-vector space $V$ with an automorphism $T$, let $M$ be a logarithm of $T$, i.e. an automorphism of $V$ such that $exp(-2\pi iM)=T$. Set
\begin{equation*}
    \rh{V}{T}:=\bigl ( V\drl u\drr,d+M\frac{d}{du} \bigr ).
\end{equation*}
For any $c \in \Cc$, consider the following $\Cc\drl u \drr$-vector space with connection
\begin{equation*}
    \widehat{\Es}^{-\frac{c}{u}}:=(\Cc \drl u \drr, d+c\frac{du}{u^2}).
\end{equation*}
Let $X$ be a smooth quasi-projective algebraic variety over $\C$ and let $W:X\rightarrow \aff{1}{\C}$ be a regular function.
For any $c \in \C$, label $\Phi_{W-c}(\C_X)$ the sheaf of vanishing cycles over the fiber $W^{-1}(c)$ and let $T$ denote the monodromy operator.

There is an equivalence of $\Z$-graded bundles with connections on $Spf(\C \drl u \drr)$
\begin{equation*}
    H^{\bullet}_{Zar}\bigl ( X, (\Omega_X^{\bullet}\drl u \drr, -dW +ud_{dR}),\nabla_{u}\bigr )\simeq \bigoplus_{c\in \C}\widehat{\Es}^{-\frac{c}{u}}\otimes_{\C \drl u \drr} \rh{\Hu^{\bullet-1}_{an}(W^{-1}(c),\Phi_{W-c}(\C_X))}{T},
\end{equation*}
where $\nabla_u=\frac{d}{du}+\frac{W}{u^2}$.
\end{thm}
\begin{rmk}
The theorem above is stated in greater generality in \cite{sa12} and \cite{sasa14}. However, we will need it only with this extent of generality.
\end{rmk}
\begin{rmk}
Notice that the sum on the right hand side is finite, as $\Phi_{W-c}(\C_X)\simeq 0$ unless $c$ is a critical value of $W$.
\end{rmk}
As an immediate consequence of the two theorems, we get an identification of periodic cyclic homology of $\mfc{X}{W}$ with vanishing cohomology (see \cite[Theorem 1.1]{ef18}) when $X$ is smooth and the critical values are contained in the fiber over $0$. In order to prove Theorem \ref{theorem Efimov}, Efimov first proves it in the affine case and then globalizes the result using a sheafification procedure, following ideas of Keller (see \cite{ke98}). We shall try to move some step towards the computation of the Hochschild and periodic cyclic homologies of $\mfc{X}{W}$ in the non smooth case, following closely Efimov's approach. We shall restrict ourselves to the affine case.

Let $A$ be a $\C$-algebra of finite type and let $g$ be a regular function on $X=Spec(A)$. We fix a presentation $A=\frac{\C[x_1,\dots,x_n]}{(f_1,\dots,f_m)}$. Set $Q=\C[x_1,\dots,x_n]$, $Y=Spec(Q)$ and $X_0=Spec(A/(g))$. Moreover, we assume that $(f_1,\dots,f_m,g)\in \C[x_1,\dots,x_n]$ is a regular sequence.

With the same notation of the previous sections, we consider the LG model $(Y,(f_1,\dots,f_m,g))$ over $(Y,\Oc_Y^{m+1})$. The global section of $\Oc_{\proj{m}{Q}}(1)$ corresponding to $(f_1,\dots,f_m,g)$ is $f_1\cdot T_1+\dots +f_m\cdot T_m+g\cdot T_{m+1}$. Similarly, the global section of $\Oc_{\proj{m-1}{Q}}(1)$ corresponding to $(f_1,\dots,f_m)$ is $f_1\cdot T_1+\dots +f_m \cdot T_m$. Set $Z^+=\Vu(f_1\cdot T_1+\dots +f_m\cdot T_m+g\cdot T_{m+1})$ and $Z=\Vu(f_1\cdot T_1+\dots +f_m\cdot T_m)$. Then diagram \ref{main diagram} in this setup becomes
\begin{equation*}
 \begindc{\commdiag}[16]
    \obj(0,50)[a]{$\Uc:=Z^+-Z$}
    \obj(60,50)[b]{$\aff{m}{Q}\simeq \aff{n+m}{\C}$}
    \obj(-60,30)[1]{$\proj{m}{X_0}$}
    \obj(0,30)[2]{$Z^+$}
    \obj(60,30)[3]{$\proj{m}{Q}$}
    \obj(-60,10)[4]{$X_0$}
    \obj(30,0)[5]{$Y$}
    \obj(-60,-10)[6]{$X$}
    \obj(-60,-30)[7]{$\proj{m-1}{X}$}
    \obj(0,-30)[8]{$Z$}
    \obj(60,-30)[9]{$\proj{m-1}{Q}=\Vu(T_{m+1}).$}
    \mor{1}{2}{$k$}
    \mor{2}{3}{$$}
    \mor{1}{4}{$q$}[\atright,\solidarrow]
    \mor{4}{6}{$i$}[\atright,\solidarrow]
    \mor(-55,10)(3,10){$$}[\atright, \solidline]
    \cmor((2,10)(15,10)(28,10)(30,8)(30,4)) \pdown(-30,10){$$}
    \mor(-55,-10)(3,-10){$$}[\atright, \solidline]
    \cmor((2,-10)(15,-10)(28,-10)(30,-8)(30,-4)) \pup(-30,-10){$$}
    \mor{7}{6}{$p$}
    \mor{7}{8}{$j$}
    \mor{8}{9}{$$}
    \mor{8}{2}{$r$}
    \mor{9}{3}{$$}
    \mor{3}{5}{$$}
    \mor{9}{5}{$$}
    \mor{a}{2}{$$}
    \mor{b}{3}{$$}
    \mor{a}{b}{$$}
 \enddc
\end{equation*}
Let $W:\aff{m}{Q}\simeq \aff{n+m}{\C}\rightarrow \aff{1}{\C}$ be the function $f_1\cdot x_{n+1}+\dots+f_m\cdot x_{n+m}+g \in \C[x_1,\dots,x_{n+m}]$, where $x_{n+k}=T_k/T_{m+1}$ for $k\in \{1,\dots,m\}$. In particular, $\Uc = \Vu(W)\subseteq \aff{n+m}{\C}$.
\begin{prop}
There is an equivalence of mixed complexes with $u$-connections
\begin{equation*}
   \bigoplus_{c\in \C} \hoch{\mfc{X}{g-c}}{b,B,\nabla^{GM}_u} \simeq \Gamma \bigl ( \aff{n+m}{\C},(\Omega_{\aff{n+m}{\C}}^{\bullet},-dW,d_{dR},\nabla_u^{dR}) \bigr ).
\end{equation*}
\end{prop}
\begin{proof}
The proof follows closely the proof of \cite[\S 4]{ef18}. 

Notice that the assignment $c\rightsquigarrow g-c$ is compatible with the assignment $g\rightarrow W$, i.e. $g-c\rightsquigarrow W-c$ for any complex number $c\in \C$. Then, by Corollary \ref{main equivalence}, we have equivalences
\begin{equation*}
    \mfc{X}{g-c}\simeq \singc{X,g-c}\simeq \sing{\aff{n+m}{\C},W-c}, \hspace{0.5cm} c\in \C.
\end{equation*}

As Hochshild homology is Morita invariant, we get that 
\begin{equation*}
    \textup{Hoch}(\mfc{X}{g-c})\simeq \textup{Hoch}(\sing{\aff{n+m}{\C},W-c})
\end{equation*}
for any $c\in \C$. Moreover, this quasi-isomorphism is strictly compatible with $u$-connections (see \cite[Proposition 3.7]{ef18}).

By \cite[\S4.10]{pp12}, we get that
\begin{equation*}
    \textup{Hoch}^{II}(\sing{\aff{n+m}{\C},W})\simeq \bigoplus_{c\in \C}\textup{Hoch}(\sing{\aff{n+m}{\C},W-c}).
\end{equation*}
Moreover, this is also strictly compatible with $u$-connections: it is defined by maps
\begin{equation*}
    \textup{Hoch}(\sing{\aff{n+m}{\C},W-c})\rightarrow  \textup{Hoch}^{II}(\sing{\aff{n+m}{\C},W-c})\xrightarrow{\simeq} \textup{Hoch}^{II}(\sing{\aff{n+m}{\C},W}),
\end{equation*}
where the first map is a morphism of mixed complexes strictly compatible with $u$-connections by \cite[Proposition 3.11]{ef18}. For what concerns the second map, let $R=\C[x_1,\dots,x_{n+m}]$ and let $R_{W}$ (resp. $R_{W-c}$, $c \in \C$) denote the $\Z/2\Z$-graded cdg algebra associated to $(R,W)$ (resp. $(R,W-c)$). Let $\sqrt{c}$ be a square root of $c$. It is immediate to see that 
\begin{equation*}
    (id,\sqrt{c}) :R_{W}\rightarrow R_{W-c}
\end{equation*}
is a cdg functor which, moreover, is a pseudo-equivalence (in the sense of \cite{pp12}). Therefore, by \cite[Proposition 2.13, Proposition 3.10, Proposition 3.13]{ef18}, the second map is a quasi-isomorphism of mixed complexes compatible with $u$-connections as well.

By Proposition 3.10, Proposition 2.13, Proposition 3.13 and Proposition 3.14 (notice that here  the assumption $crit(W)\subseteq W^{-1}(0)$ is not needed) in \cite{ef18}, we get that 
\begin{equation*}
    \textup{Hoch}^{II}(\sing{\aff{n+m}{\C},W})\simeq \Gamma \bigl ( \aff{n+m}{\C},(\Omega_{\aff{n+m}{\C}}^{\bullet},-dW,d_{dR},\nabla_u^{dR}) \bigr ).
\end{equation*}
\end{proof}

In order to have an equivalence in the statement of the previous proposition with just the term $\hh{\mfc{X}{g}{,b,B,\nabla^{GM}_u}}$ on the left hand side, one needs to impose some hypothesis, similar to the requirement that the critical locus of $g$ is contained in the fiber over zero when $X$ is smooth. For this we introduce the following 
\begin{defn}\label{relative critical locus}
Let $g:X\rightarrow \aff{1}{\C}$ be as above. We say that $x\in X$ is a \textit{relative critical point} of $W$ if
\begin{equation*}
    rank(Jac(f_1,\dots,f_m))(x)=rank(Jac(f_1,\dots,f_m,g))(x).
\end{equation*}
A \textit{relative critical value} $c\in \aff{1}{\C}$ is the image of a relative critical point.
\end{defn}
\begin{rmk}
If $X$ is smooth, a relative critical point is just a critical point in the usual sense.
\end{rmk}
\begin{rmk}
It seems possible that Definition \ref{relative critical locus} is related to the notion of critical value of a regular function on a singular variety, introduced in \cite[\S B.2]{ep15}.
\end{rmk}
\begin{ex}
Let $z:X=Spec(\C[x,y,z]/(y^2-x^3))\rightarrow \aff{1}{\C}$. All points of the line $\Vu(x,y)$ are critical points (meaning that all the fibers are singular), but none of them is a relative critical point.
\end{ex}
\begin{prop}
Let $(X,g)$ be as above. Assume that $0$ is the only relative critical value. Then
\begin{equation*}
    \hoch{\mfc{A}{g}}{b,B,\nabla}\simeq \Gamma \bigl ( \aff{n+m}{\C},(\Omega_{\aff{n+m}{\C}}^{\bullet},-dW,d_{dR},\nabla_u^{dR}) \bigr )
\end{equation*}
is an equivalence of mixed complexes with $u$-connection.
\end{prop}
\begin{proof}
The hypothesis on the relative critical locus of $g$ means that the Jacobian of $W$
\begin{equation*}
    Jac(W)=\Bigl ( \frac{dg}{dx_1}+\sum_{k=1}^m\frac{df_k}{dx_1}x_{n+k},\dots,\frac{dg}{dx_n}+\sum_{k=1}^m\frac{df_k}{dx_n}x_{n+k},f_1,\dots,f_m\Bigr )
\end{equation*}
can vanish only if $f_1(x_1,\dots,x_n)=\dots= f_m(x_1,\dots,x_n)=g(x_1,\dots,x_n)=0$, in which case $W=0$. In particular,
\begin{equation*}
    \textup{Hoch}^{II}(\sing{\aff{n+m}{\C},W})\simeq \textup{Hoch}(\sing{\aff{n+m}{\C},W})
\end{equation*}
and the statement is clear from the previous proposition.
\end{proof}
As an immediate consequence of the two propositions above, we have the following
\begin{thm}\label{computation HH and HP mf coh}
With the same notation as above, there is an equivalence of $\Z/2\Z$-graded bundles with connections on $Spf(\C \drl u \drr)$
\begin{equation*}
    \bigoplus_{c\in \C} \hp{\mfc{X}{g-c},\nabla^{GM}_u}\simeq \bigl ( \Hu^{\bullet}_{Zar}(\Omega^{\bullet}_{\aff{n+m}{\C}}\drl u \drr, -dW+ud_{dR}), \nabla_u^{dR}=\frac{d}{du}+\frac{\Gamma}{u}+\frac{W}{u^2} \bigr ),
\end{equation*}
where $\Gamma_{|\Omega^p_{\aff{n+m}{\C}}}=-\frac{p}{2}\cdot id$, and an equivalence of $\Z/2\Z$-graded vector spaces
\begin{equation*}
    \bigoplus_{c\in \C}\hh{\mfc{X}{g-c}}\simeq \Hu^{\bullet}_{Zar}(\aff{n+m}{\C},(\Omega^{\bullet}_{\aff{n+m}{\C}},-dW)).
\end{equation*}
Moreover, if the relative critical locus of $g$ is contained in the fiber over $0$, then
\begin{equation*}
   \hp{\mfc{X}{g},\nabla^{GM}_u}\simeq \bigl ( \Hu^{\bullet}_{Zar}(\Omega^{\bullet}_{\aff{n+m}{\C}}\drl u \drr, -dW+ud_{dR}), \nabla_u^{dR}=\frac{d}{du}+\frac{\Gamma}{u}+\frac{W}{u^2} \bigr ),
\end{equation*}
and 
\begin{equation*}
    \hh{\mfc{X}{g}}\simeq \Hu^{\bullet}_{Zar}(\aff{n+m}{\C},(\Omega^{\bullet}_{\aff{n+m}{\C}},-dW)).
\end{equation*}
\end{thm}
This theorem, combined with Sabbah/Sabbah-Saito's result, gives us the following computation of the periodic cyclic homology of $\mfc{X}{g}$ in terms of vanishing cohomology in the affine case.
\begin{thm}\label{HP and vanishing cycles}
With the same notation as above, there is an equivalence of $\Z/2\Z$-graded vector bundles with connections of the punctured disk, 
\begin{equation*}
    \bigoplus_{c\in \C}\hp{\mfc{X}{g-c},\nabla^{GM}_u}\simeq \bigoplus_{c\in \C}\widehat{\Es}^{-\frac{c}{u}}\otimes_{\C \drl u \drr} \rh{\Hu^{\bullet-1}_{an}(W^{-1}(c),\Phi_{W-c}(\C_{\aff{n+m}{\C}}))}{T\cdot (-1)^{\bullet}}.
\end{equation*}
Moreover, if the relative critical locus of $g$ is contained in the fiber over $0$, then
\begin{equation*}
    \hp{\mfc{X}{g},\nabla^{GM}_u}\simeq  \rh{\Hu^{\bullet-1}_{an}(W^{-1}(0),\Phi_{W}(\C_{\aff{n+m}{\C}}))}{T\cdot (-1)^{\bullet}}.
\end{equation*}
\end{thm}
\begin{rmk}
The $(-1)^{\bullet}$ in front of the monodromy operator $T$ appears as we have added the term $\frac{\Gamma}{u}$ to the connection on the twisted de Rham complex. See the proof of Theorem 5.4 in \cite{ef18}.
\end{rmk}
\begin{rmk}
In order to obtain a generalization of the theorem above to the non affine case, the equivalence
\begin{equation*}
    \mfc{X}{s}\simeq \sing{\Vc,W_{|\Vc}\in \Gamma(\Vc,\Oc(1))}
\end{equation*}
suggests that a formalism of vanishing cycles over $\ag{\C}$ analogous to the one sketched in \cite{p20} is required. 

It also seems possible that there exists an algebraic computation for vanishing cycles over $\ag{\C}$ that generalizes Kontsevich's formula.

This is currently being investigated by the author.
\end{rmk}
\begin{rmk}
In \cite[\S 6]{ef18}, the author states that the following formula is expected to hold
\begin{equation*}
    \bigl ( \hp{\mfc{X}{g}},\nabla_u^{GM} \bigr )\simeq \rh{(\Hu^{\bullet-1}_c(g^{-1}(0),\Phi_{g}(\C_X))^{\vee}}{T^{\vee}\cdot (-1)^{\bullet}},
\end{equation*}
at least when the relative critical locus of $W$ is concentrated over $0$.

In light of the theorem above, this is equivalent to proving that there is a quasi-isomorphism
\begin{equation*}
    \rh{\Hu^{\bullet-1}_{an}(W^{-1}(0),\Phi_{W}(\C_{\aff{n+m}{\C}}))}{T\cdot (-1)^{\bullet}}\simeq \rh{(\Hu^{\bullet-1}_c(g^{-1}(0),\Phi_{g}(\C_X))^{\vee}}{T^{\vee}\cdot (-1)^{\bullet}}.
\end{equation*}
\end{rmk}

\end{document}